\documentclass[11pt]{amsart} 
\usepackage{amssymb,tensor,enumitem,verbatim}
\usepackage[bookmarksdepth=3]{hyperref}
\usepackage[margin=1in]{geometry}
\usepackage[all,cmtip]{xy}
\usepackage[sort,noadjust]{cite}
\usepackage[color=yellow]{todonotes}

\newcommand{\subgrp}[1]{\langle #1 \rangle}
\newcommand{\set}[1]{\left\{ #1 \right\}}
\newcommand{\abs}[1]{\left| #1 \right|}
\newcommand{\bs}[1]{\boldsymbol{#1}}

\newcommand{\wt}[1]{\widetilde{ #1}}
\newcommand{\ol}[1]{\overline{#1}}

\newcommand{\gotimes}{\tensor[^g]{\otimes}{}}

\DeclareMathOperator{\ann}{ann}

\newcommand{\Cbul}{C^\bullet}

\DeclareMathOperator{\cx}{cx}

\DeclareMathOperator{\Dist}{Dist}

\DeclareMathOperator{\Ext}{Ext}

\DeclareMathOperator{\gldim}{gldim}
\DeclareMathOperator{\opH}{H}
\newcommand{\Hbul}{\opH^\bullet}
\DeclareMathOperator{\Hom}{Hom}

\DeclareMathOperator{\Lie}{Lie}

\DeclareMathOperator{\Proj}{Proj}
\newcommand{\NProj}{N\text{-}\Proj}
\DeclareMathOperator{\projdim}{projdim}
\DeclareMathOperator{\MaxSpec}{MaxSpec}
\DeclareMathOperator{\rank}{rank}
\DeclareMathOperator{\res}{res}

\DeclareMathOperator{\Spec}{Spec}
\DeclareMathOperator{\supp}{supp}
\newcommand{\pisupp}{\pi\text{-}\supp}

\newcommand{\ev}{\textup{ev}}
\newcommand{\odd}{\textup{odd}}
\newcommand{\ve}{\varepsilon}
\newcommand{\zero}{\ol{0}}
\newcommand{\one}{\ol{1}}

\newcommand{\calN}{\mathcal{N}}

\newcommand{\Chi}{\mathcal{X}}

\newcommand{\bsa}{\bs{A}}
\newcommand{\bse}{\bs{e}}

\newcommand{\bsc}{\bs{c}}
\newcommand{\bsl}{\bs{\Lambda}}

\newcommand{\bsV}{\bs{V}}
\newcommand{\bsvr}{\bsV\!_r}

\newcommand{\C}{\mathbb{C}}

\newcommand{\G}{\mathbb{G}}
\newcommand{\M}{\mathbb{M}}
\newcommand{\N}{\mathbb{N}}
\renewcommand{\P}{\mathbb{P}}
\newcommand{\Z}{\mathbb{Z}}

\newcommand{\Ga}{\G_a}
\newcommand{\Gam}{\Ga^-}
\newcommand{\Gar}{\G_{a(r)}}
\newcommand{\Gas}{\G_{a(s)}}

\newcommand{\GLn}{GL_n}
\newcommand{\GLnr}{GL_{n(r)}}
\newcommand{\GLmn}{GL_{m|n}}
\newcommand{\GLmnr}{GL_{m|n(r)}}

\newcommand{\Mr}{\M_r}
\newcommand{\Mrfeta}{\M_{r,f,\eta}}
\newcommand{\Mrseta}{\M_{r,s,\eta}}
\newcommand{\Msfeta}{\M_{s,f,\eta}}
\newcommand{\Mrs}{\M_{r,s}}
\renewcommand{\Pr}{\P_r}
\newcommand{\Pone}{\P_1}
\newcommand{\Vrg}{V_r(G)}
\newcommand{\Vrs}{V_{r,s}}
\newcommand{\Vrfeta}{V_{r,f,\eta}}

\newcommand{\g}{\mathfrak{g}}
\newcommand{\gl}{\mathfrak{gl}}
\newcommand{\glmn}{\gl(m|n)}
\newcommand{\gzero}{\g_{\zero}}
\newcommand{\gone}{\g_{\one}}

\newcommand{\wtg}{\wt{\g}}

\newcommand{\fb}{\mathfrak{b}}
\newcommand{\fp}{\mathfrak{p}}

\newcommand{\fs}{\mathfrak{s}}

\newcommand{\ft}{\mathfrak{t}}

\newcommand{\Aone}{A_{\one}}
\newcommand{\Lzero}{L_{\zero}}
\newcommand{\Lone}{L_{\one}}
\newcommand{\Vone}{V_{\ol{1}}}

\newcommand{\Vzero}{V_{\ol{0}}}

\newcommand{\CAlg}{\mathsf{CAlg}}
\newcommand{\CSAlg}{\mathsf{CSAlg}}

\newcommand{\sHopf}{\mathsf{sHopf}}
\newcommand{\Set}{\mathsf{Set}}
\newcommand{\Grp}{\mathsf{Grp}}
\newcommand{\sGrp}{\mathsf{sGrp}}

\newcommand{\sMod}{\mathsf{sMod}}

\newcommand{\fg}{\mathfrak{g}}
\newcommand{\ff}{\mathfrak{f}}
\newcommand{\fe}{\mathfrak{e}}
\newcommand{\fa}{\mathfrak{a}}
\newcommand{\fh}{\mathfrak{h}}
\newcommand{\osp}{\mathfrak{osp}}

\newcommand{\FF}{\mathcal{F}}
\newcommand{\VV}{\mathcal{V}}
\newcommand{\KK}{\mathcal{K}}

\newcommand{\0}{\zero}
\renewcommand{\1}{\one}

\DeclareMathOperator{\atyp}{atyp}
\DeclareMathOperator{\defect}{def}
\DeclareMathOperator{\Norm}{Norm}
\DeclareMathOperator{\sdim}{sdim}
\DeclareMathOperator{\Stab}{Stab}
\DeclareMathOperator{\strace}{str}

\numberwithin{equation}{subsection}

\newtheorem{theorem}{Theorem}[subsection]
\newtheorem{proposition}[theorem]{Proposition}

\theoremstyle{definition}
\newtheorem{definition}[theorem]{Definition}

\newtheorem{question}[theorem]{Question}

\newtheorem{problem}[theorem]{Problem}
\newtheorem{remark}[theorem]{Remark}

\title[Support theories for Lie super\-algebras and finite super\-group schemes]{A survey of support theories for Lie super\-algebras and finite super\-group schemes}

\author{Christopher M.\ Drupieski}
\address{Department of Mathematical Sciences,
		DePaul University,
		Chicago, IL 60614, USA}
\email{c.drupieski@depaul.edu}

\thanks{The first author was supported in part by Simons Collaboration Grant for Mathematicians No.\ 4269055.}

\author{Jonathan R.\  Kujawa}
\address{Department of Mathematics \\
		University of Oklahoma \\
		Norman, OK 73019, USA}
\email{kujawa@math.ou.edu}

\thanks{The second author was supported in part by Simons Collaboration Grant for Mathematicians No.\ 525043.}

\date{\today}



\allowdisplaybreaks

\begin{document}

\begin{abstract}
We survey the current state of various support variety theories for Lie super\-algebras and finite supergroup schemes. We pay particular attention to the theory in characteristic zero developed by Boe, Kujawa, and Nakano using relative Lie superalgebra cohomology, and to the theory developed in positive characteristic in our previous work.
\end{abstract}

\maketitle

\setcounter{tocdepth}{1}
\tableofcontents

\section{Introduction}

\subsection{Overview}

Let $A$ be a Hopf algebra or similar algebraic structure over an algebraically closed field $k$. Support varieties are tools that help relate the representation theory of $A$ to the geometry encoded by the spectrum $\abs{A}$ of the cohomology ring $\Hbul(A,k) = \Ext_A^\bullet(k,k)$. The enveloping algebra $U(L)$ of a finite-dimensional Lie algebra $L$ has finite global dimension $d = \dim(L)$, so the Lie algebra cohomology ring $\Hbul(L,k) = \Hbul(U(L),k)$ vanishes in all large degrees, and hence its spectrum does not provide an interesting geometry. On the other hand, in positive characteristic the cohomology of finite-dimensional \emph{restricted} Lie algebras---or more generally, of finite group schemes---frequently gives rise to rich algebro-geometric structures. Support variety theories in these contexts have been developed and explored by Friedlander and Parshall \cite{Friedlander:1986b}, Suslin, Friedlander, and Bendel \cite{Suslin:1997a}, and Friedlander and Pevtsova \cite{Friedlander:2007}, and many others.

The aim of this paper is to give an overview of what is currently known about cohomological support varieties for Lie \emph{super}algebras and finite \emph{super}group schemes. For a finite-dimensional Lie superalgebra $L$, the enveloping super\-algebra $U(L)$ may or may not have finite global dimension, depending on whether or not $L$ contains any odd elements such that $[x,x] = 0$. In characteristic zero, even if a Lie super\-algebra $L$ has infinite global dimension, the ordinary Lie super\-algebra cohomology ring $\Hbul(L,k)$ may still be finite, and hence may not give rise to any interesting ambient geometry; see for example the calculations in \cite[\S2.6]{Fuks:1986}. But replacing the ordinary Lie superalgebra cohomology ring with the relative cohomology ring for the pair $(L,\Lzero)$, one gets a rich support variety theory in many cases. In positive characteristic, the Lie super\-algebra cohomology ring $\Hbul(L,k)$ gives rise to a nontrivial geometry whenever the enveloping superalgebra $U(L)$ is not of finite global dimension. In general, support varieties for finite-dimensional restricted Lie super\-algebras (and more generally, for finite super\-group schemes) are not completely understood, with the fullest picture available only for $p$-nilpotent restricted Lie superalgebras (and more generally, for unipotent finite supergroup schemes). What is known at present parallels the classical, non-super theory in many ways, but with interesting twists arising from new super phenomena.

The paper is organized as follows: We begin in Section \ref{sec:support-and-cohomology} by recalling some of the basic definitions and properties for support varieties and Lie superalgebra cohomology. Part \ref{part:char0} of the paper is devoted to what is known in characteristic zero. In Section \ref{sec:support-via-relative} we summarize the support variety theory for Lie super\-algebras developed in a series of papers by Boe, Kujawa, Nakano, and collaborators using relative Lie super\-algebra cohomology. In Section \ref{sec:supergroup-char-0} we give a brief account of support varieties for finite super\-group schemes in characteristic zero. Part \ref{part:char-p} of the paper is devoted to what is known in (odd) positive characteristic. In Section \ref{sec:support-ordinary} we describe the support variety theory one gets for a finite-dimensional Lie superalgebra $L$ from considering its ordinary Lie superalgebra cohomology ring $\Hbul(L,k)$. Then, to aid the reader in comparing the super theory with its classical (non-super) counterpart, in Section \ref{sec:non-super} we summarize the classical results for restricted Lie algebras and infinitesimal group schemes (as developed by Friedlander and Parshall and by Suslin, Friedlander, and Bendel), before discussing in Section \ref{sec:infinitesimal-supergroups} our work to date related to restricted Lie superalgebras and infinitesimal supergroup schemes. As a complement to our own work on infinitesimal supergroups, we also recommend to the reader the recent survey article by Benson, Iyengar, Krause, and Pevtsova \cite{Benson:2022}.    

\subsection{Conventions}

Throughout, $k$ will denote a ground field of characteristic $p \neq 2$. For simplicity, we will assume that $k$ is algebraically closed, though it may be possible to formulate some results without this general assumption. Except when specified, all vector spaces will be $k$-vector spaces, all algebras will be (associative, unital) $k$-algebras, and all unadorned tensor products will be tensor products over $k$. Given a $k$-vector space $V$, let $V^* = \Hom_k(V,k)$ be its $k$-linear dual.

Set $\Z_2 = \Z/2\Z = \{ \zero,\one \}$. Following the literature, we use the prefix `super' to indicate that an object is $\Z_2$-graded, and we assume that the reader is generally familiar with the standard conventions of `super\-mathematics'. In particular, we denote the decomposition of a vector superspace into its $\Z_2$-homogeneous components by $V = \Vzero \oplus \Vone$, calling $\Vzero$ and $\Vone$ the even and odd subspaces of $V$, respectively, and writing $\ol{v} \in \Z_2$ to denote the superdegree of a homogeneous element $v \in \Vzero \cup \Vone$. When written without additional adornment, we consider the ground field to be a superspace concentrated in even super\-degree. Whenever we state a formula in which homogeneous degrees of elements are specified, we mean that the formula is true as written for homogeneous elements and that it extends by linearity to non-homogeneous elements.

We use the adjective \emph{graded} to indicate that an object admits an additional $\Z$-grading that is compatible with its underlying structure. Thus a \emph{graded superspace} is a $(\Z \times \Z_2)$-graded vector space, a \emph{graded super\-algebra} is a $(\Z \times \Z_2)$-graded algebra, etc. Given a graded superspace $V$ and a homogeneous element $v \in V$ of bidegree $(s,t) \in \Z \times \Z_2$, we write $\deg(v) = s$ and $\ol{v} = t$ for the $\Z$-degree and the superdegree of $v$, respectively. If $A$ is a graded super\-algebra, we say that $A$ is \emph{graded-commutative} provided that for all homogeneous elements $a,b \in A$, one has
	\begin{equation} \label{eq:graded-commutative}
	ab = (-1)^{\deg(a) \cdot \deg(b) + \ol{a} \cdot \ol{b}}ba.
	\end{equation}
The notion of \emph{graded-cocommutativity} for a graded super-coalgebra is defined similarly. If $A$ and $B$ are graded super\-algebras, then $A \otimes B$ is a graded super\-algebra, with $\deg(a \otimes b) = \deg(a) + \deg(b)$, and with the product defined by $(a \otimes b) (c \otimes d) = (-1)^{\deg(b) \cdot \deg(c) + \ol{b}\cdot \ol{c}} ac \otimes bd$; we denote this graded super\-algebra by $A \gotimes B$, and call it the \emph{graded tensor product of super\-algebras}.

Let $\N = \set{0,1,2,3,\ldots}$ be the set of non-negative integers.

\section{Background on support varieties and Lie superalgebra cohomology} \label{sec:support-and-cohomology}

\subsection{Basic notions for support varieties} \label{subsec:basic-notions-support}

In this section we recall some of the basic definitions concerning (cohomological) support varieties for supermodules over a Hopf super\-algebra. For further details, the reader can consult \cite[\S2]{Drupieski:2019a}.

Let $k$ be a (for simplicity, algebraically closed) field of characteristic $p \neq 2$, and let $A$ be a Hopf super\-algebra over $k$. The cohomology ring $\Hbul(A,k) = \Ext_A^\bullet(k,k)$ is naturally a graded super\-algebra: the $\Z$-component of the grading is the cohomological grading, and the $\Z_2$-component is the `internal' grading coming from the super\-algebra structure on $A$. The cup product makes $\Hbul(A,k)$ into a graded-commutative super\-algebra in the sense of \eqref{eq:graded-commutative}. In particular, the subspaces $\opH^{\ev}(A,k)_{\one}$ and $\opH^{\odd}(A,k)_{\zero}$ of $\Hbul(A,k)$ consist of nilpotent elements. Then modulo its nilradical, $\Hbul(A,k)$ is a commutative ring in the ordinary sense, with the super\-degree of a homogeneous element equal simply to the reduction modulo $2$ of its $\Z$-degree.

Given $A$-supermodules $M$ and $N$, the cohomology ring $\Hbul(A,k)$ acts on $\Ext_A^\bullet(M,N)$ via the cup product. In all of the examples we will consider, the following properties will hold:
	\begin{enumerate}[wide=\parindent,label={(FG\arabic*)}]
	\item \label{FG-ring} $\Hbul(A,k)$ is a finitely-generated $k$-algebra
	\item \label{FG-modules} If $M$ and $N$ are finite-dimensional, then $\Ext_A^\bullet(M,N)$ is a finite $\Hbul(A,k)$-module.
	\end{enumerate}
Let $I_A(M)$ be the annihilator ideal for the left cup product action of $\Hbul(A,k)$ on $\Ext_A^\bullet(M,M)$:
	\[
	I_A(M) = \ann_{\Hbul(A,k)} \Ext_A^\bullet(M,M).
	\]
Equivalently, $I_A(M)$ is the kernel of the algebra map $\Phi_M : \Ext_A^\bullet(k,k) \to \Ext_A^\bullet(M,M)$ induced by the functor $- \otimes M$. Now the cohomological spectrum of $A$ is
	\[
	\abs{A} = \MaxSpec \Hbul(A,k),
	\]
the maximal ideal spectrum of $\Hbul(A,k)$, and the support variety of $M$, denoted $\abs{A}_M$, is the Zariski closed subset of $\abs{A}$ defined by $I_A(M)$:
	\[
	\abs{A}_M = \MaxSpec \big( \Hbul(A,k)/I_A(M) \big).
	\]
Support varieties for supermodules over a Hopf super\-algebra satisfy a list of standard properties; see for example \cite[\S2.3]{Drupieski:2019a} or \cite[\S8.3]{Witherspoon:2019}. In particular, if $A$ is a finite-dimensional (hence self-injective) Hopf super\-algebra, then under the `finite-generation' hypotheses \ref{FG-ring} and \ref{FG-modules}, one gets for a finite-dimensional $A$-supermodule $M$ that the geometric dimension of the support variety $\abs{A}_M$ is equal to $\cx_A(M)$, the complexity of $M$ as an $A$-supermodule, and the support variety $\abs{A}_M$ is zero (i.e., is equal to a single point) if and only if $M$ is projective; see \cite[Proposition 2.3.13]{Drupieski:2019a}.

In some contexts we may instead consider support varieties defined via either the prime ideal spectrum $\Spec \Hbul(A,k)$ or the projective spectrum $\Proj \Hbul(A,k)$. When confusion is unlikely, we may write $\abs{A} = \Spec \Hbul(A,k)$ and $\abs{A}_M = \Spec \big( \Hbul(A,k)/I_A(M) \big)$.

In Part \ref{part:char0} we will consider support varieties that are defined via relative cohomology groups. Aside from the switch to relative cohomology, the basic definitions of the cohomological spectrum and of support varieties are the same. Support varieties defined via relative cohomology may not satisfy the full list of `standard' properties as those defined via non-relative cohomology.

\subsection{Generalities of Lie super\-algebra cohomology} \label{subsec:generalities-cohomology}

Let $L$ be a Lie super\-algebra over the field $k$, and let $U(L)$ be the universal enveloping super\-algebra of $L$. By construction, an $L$-super\-module is the same thing as a $U(L)$-supermodule. The field $k$ is an $L$-supermodule via the augmentation map $\ve: U(L) \to k$. Given $L$-supermodules $M$ and $N$, set
	\[
	\Ext_L^\bullet(M,N) = \Ext_{U(L)}^\bullet(M,N) \quad \text{and} \quad \Hbul(L,M) = \Ext_L^\bullet(k,M).
	\]
By definition, the differentials in a $U(L)$-supermodule resolution are required to be even linear maps, i.e., we only consider resolutions in $(\sMod)_{\ev}$, the `underlying even subcategory' of the category $\sMod_{U(L)}$ of all $U(L)$-super\-modules. Then $\Ext_L^\bullet(M,N)$ can be computed as the cohomology of the complex $\Hom_{U(L)}(P_\bullet,N)$, or as the cohomology of the complex $\Hom_{U(L)}(M,Q^\bullet)$, for any projective $U(L)$-supermodule resolution $P_\bullet \to M$, or any injective $U(L)$-supermodule resolution $N \to Q^\bullet$.

An explicit $U(L)$-free resolution of the ground field is provided by the Koszul resolution, denoted $Y(L)$. As a left $U(L)$-supermodule, $Y(L) = U(L) \otimes \bsa(L)$, where $\bsa(L) = \bigoplus_{n \in \N} \bsa^n(L)$ denotes the alternating power super\-algebra on $L$, as defined for example in \cite[\S2.3.7]{Drupieski:2016}. As a graded super\-algebra, $\bsa(L) \cong \Lambda(\Lzero) \gotimes \Gamma(\Lone)$. Here $\Lambda(\Lzero) = \bigoplus_{n \in \N} \Lambda^n(\Lzero)$ is the ordinary exterior algebra on $\Lzero$ (considered as a graded super\-algebra concentrated in superdegree $\zero$), and $\Gamma(\Lone) = \bigoplus_{n \in \N} \Gamma^n(\Lone)$ is the ordinary divided power algebra on $\Lone$ (considered as a graded super\-algebra with $\Gamma^n(\Lone)$ concentrated in superdegree $\ol{n}$). One can define a graded super-bialgebra structure on $Y(L)$ such that the Koszul differential makes $Y(L)$ into a differential graded super-bialgebra, and such that the coproduct $Y(L) \to Y(L) \gotimes Y(L)$ restricts to the usual coproducts on $U(L)$, $\Lambda(\Lzero)$, and $\Gamma(\Lone)$. For more details, see \cite[\S3.1]{Drupieski:2013c}. In particular, an explicit formula for the Koszul differential $d: Y_n(L) \to Y_{n-1}(L)$ is given in \cite[Remark 3.1.4]{Drupieski:2013c}.

\begin{remark} \label{rem:subcomplex-bar}
By definition, $\bsa^n(L)$ is a subspace of $L^{\otimes n}$. Then via the canonical maps $\bsa^n(L) \subseteq L^{\otimes n} \subseteq U(L)^{\otimes n}$, one can show that $Y(L)$ is a subcomplex of the left bar complex for $U(L)$.
\end{remark}

In the special case that $k$ is a field of characteristic zero, one has $\gamma_a(y) = \frac{1}{a!} [\gamma_1(y)]^a$ in $\Gamma(\Lone)$, and $\Gamma(\Lone)$ is isomorphic to the symmetric algebra $S(\Lone)$. Thus in characteristic zero, the underlying space for the Kosuzl resolution can be realized as $Y(L) = U(L) \otimes \bsl(L)$, where $\bsl(L)$ denotes the superexterior algebra on the superspace $L$, as discussed for example in \cite[\S2.3.5]{Drupieski:2016}.

Now given an $L$-supermodule $M$, $\Hbul(L,M)$ can be computed as the cohomology of the cochain complex $\Cbul(L,M) := \Hom_{U(L)}(Y_\bullet(L),M)$. Denote the differential on $\Cbul(L,M)$ by $\partial$. In the case $M = k$, the coproduct on $Y(L)$ induces an algebra structure on $\Cbul(L,k)$, which in turn induces the cup product on the cohomology ring $\Hbul(L,k)$. Explicitly, $\Cbul(L,k)$ is isomorphic as a graded super\-algebra to the superexterior algebra $\bsl(L^*)$. The graded super\-algebra isomorphism
	\[
	\Cbul(L,k) = \Hom_{U(L)}(U(L) \otimes \bsa(L),k) \cong \Hom_k(\bsa(L),k) \cong \bsl(L^*)
	\]
corresponds to the duality of strict polynomial superfunctors $\bsa^\# \cong \bsl$ discussed in \cite[\S2.6]{Drupieski:2016}. As a graded super\-algebra, $\bsl(L^*) \cong \Lambda(\Lzero^*) \gotimes S(\Lone^*)$, where $S(\Lone^*) = \bigoplus_{n \in \N} S^n(\Lone^*)$ denotes the ordinary symmetric algebra on $\Lone^*$ (considered as a graded super\-algebra with $S^n(\Lone^*)$ concentrated in superdegree $\ol{n}$). The differential $\partial$ makes $\Cbul(L,k)$ into a differential graded super\-algebra. 

The left adjoint action of $L$ on itself extends by superderivations to an action of $L$ on $\bsa(L)$. Considering $\bsa^n(L)$ as a subspace of $L^{\otimes n}$, this is just the restriction of the usual adjoint action on tensor space. On a monomial $\gamma_a(y) \in \Gamma(\Lone)$, the adjoint action of $u \in L$ is given by $u.\gamma_a(y) = s([u,y]) \gamma_{a-1}(y)$, where $s: L \to \bsa^1(L)$ is the canonical map. Now let $\fs$ be a Lie super\-algebra over $k$, and let $\ft$ be a Lie sub-super\-algebra of $\fs$. The adjoint action of $\fs$ on $\bsa(\fs)$ descends to an action of $\ft$ on $\bsa(\fs/\ft)$. Also, $U(\fs)$ becomes a right $\ft$-supermodule via right multiplication. Set
	\[
	Y(\fs,\ft) = U(\fs) \otimes_{U(\ft)} \bsa(\fs/\ft).
	\]
The Koszul differential on $Y(\fs)$ induces a differential on $Y(\fs,\ft)$, which together with the induced augmentation map $U(\fs) \otimes_{U(\ft)} \bsa(\fs/\ft) \to U(\fs) \otimes_{U(\ft)} k \to k$ makes $Y(\fs,\ft)$ into a resolution of the field $k$. We call this complex the \emph{relative Koszul resolution for the pair $(\fs,\ft)$}.

The relative Lie super\-algebra cohomology of the pair $(\fs,\ft)$ with coefficients in an $\fs$-supermodule $M$, denoted $\Hbul(\fs,\ft;M)$, is by definition the cohomology of the cochain complex $\Hom_{U(\fs)}(Y(\fs,\ft),M)$. If $k$ is a field of characteristic zero (e.g., if $k = \C$), and if $\fs$ is finitely-semisimple under the adjoint action of $\ft$, then $Y(\fs,\ft)$ is a $(U(\fs),U(\ft))$-projective resolution of the field $k$, and hence $\Hbul(\fs,\ft;M)$ is equal to the relative $\Ext$-group $\Ext_{(U(\fs),U(\ft))}^\bullet(k,M)$ as defined for example by Hochschild \cite{Hochschild:1956}. For more details on relative homological algebra, see for example \cite[Chapter III]{Kumar:2002}.

\hypersetup{bookmarksdepth=-1}
\part{Characteristic zero} \label{part:char0}
\hypersetup{bookmarksdepth=3}

In this part, let $k = \C$.

\section{Support varieties via relative Lie superalgebra cohomology} \label{sec:support-via-relative}

\subsection{Relative cohomology and support varieties}

A finite-dimensional complex Lie super\-algebra $\fg = \fg_{\zero} \oplus \fg_{\one}$ is called \emph{classical} if there is a connected reductive algebraic group $G_{\zero}$ such that $\Lie(G_{\zero}) = \fg_{\zero}$ and if there is an action of $G_{\zero}$ on $\fg_{\one}$ that differentiates to the adjoint action of $\fg_{\zero}$ on $\fg_{\one}$.  In this section we will assume $\fg$ is classical. Our running example will be $\fg = \glmn$. Recall that
	\begin{equation}\label{E:glmndef}
	\glmn = \left\{ \left(	\begin{array}{c|c}
								A & B \\
								\hline
								C & D
								\end{array} \right) \right\},
	\end{equation}
where $A$ is a $m \times m$ matrix, $B$ is a $m \times n$ matrix, $C$ is a $n \times m$ matrix, and $D$ is a $n \times n$ matrix. The $\Z_{2}$-grading is defined so that $\fg_{\0}$ consists of those matrices with $B=0$ and $C=0$, $\fg_{\1}$ is those matrices with $A=0$ and $D=0$, and the Lie bracket is the supercommutator.  Note that $G_{\zero} = GL_m \times GL_n$ acts on $\fg_{\1}$ by matrix conjugation, and $\glmn$ is classical.

Let $\FF = \FF (\fg , \fg_{\zero})$ be the category of finite-dimensional $\fg$-supermodules that are completely reducible when restricted to $\fg_{\zero}$.  It is common that $\fg_{\0}$ is a semisimple Lie algebra, in which case $\FF$ is the category of all finite-dimensional $\fg$-supermodules.  One can use the relative cohomology introduced in Section~\ref{subsec:generalities-cohomology} to compute cohomology in $\FF$.  Namely, given supermodules $M$ and $N$ in $\FF$, \cite[Theorem 2.5]{Boe:2010} states that for all $d \geq 0$, 
	\[
	\Ext^{d}_{\FF}(M, N) \cong \opH^{d}(\fg , \fg_{\zero}; M^{*} \otimes N)
	\]
as superspaces.  In the particular case when $M=N = k$, it is easy to see that the differentials in the relative Koszul complex are identically zero and hence the cohomology ring is a ring of invariants:
	\[
	\Ext_{\FF}^{\bullet} (k,k) \cong \Hbul (\fg , \fg_{\zero}; k) \cong S^{\bullet}(\fg^{*}_{\1})^{\fg_{\zero}}=S^{\bullet}(\fg^{*}_{\1})^{G_{\zero}}.
	\]
Since $\fg$ is classical, $\Hbul (\fg , \fg_{\zero}; k)$ is a finitely-generated algebra and $\Ext^{\bullet}_{\FF}(M,M)$ is a finitely-generated $\Hbul (\fg , \fg_{\zero}; k)$-module.  Thus conditions \ref{FG-ring} and \ref{FG-modules} are satisfied.

Using notions from invariant theory allows us to compute $\Hbul (\fg , \fg_{\zero}; k)$ in many cases. If there exists an element $x_{0} \in \fg_{\one}$ such that the orbit $G_{\zero} \cdot x_{0}$ is both closed and has maximal dimension among $G_{\zero}$-orbits, then the action of $G_{\zero}$ on $\fg_{\one}$ is said to be \emph{stable}.  Set
	\[
	\Stab_{G_{\zero}}(x_{0}) = \left\{g \in G_{\zero} : g \cdot x_{0}=x_{0} \right\}.
	\]
When the action of $G_{\zero}$ on $\fg_{\1}$ is stable one can define the detecting subalgebra $\ff = \ff_{\zero} \oplus \ff_{\one}$ via
	\[
	\ff_{\one } = \fg_{\one}^{\Stab_{G_{\zero}}(x_{0})} \quad \text{ and } \quad \ff_{\zero} = [\ff_{\one}, \ff_{\one }].
	\]
For short, we say that $\fg$ is \emph{stable} if the action of $G_{\zero}$ on $\fg_{\one}$ is stable. Set
	\[
	N = \Norm_{G_{\zero}}(\ff_{\1})=\left\{g \in G_{\zero} : g \cdot \ff_{\one} \subseteq \ff_{\one}  \right\}.
	\]

If there exists a vector $v \in \fg_{\one}$ such that $\dim_{k} \left\{ x \in \fg_{\one} : \fg_{\zero}.x \subseteq \fg_{\zero}.v \right\} = \dim S^{\bullet}(\fg_{\one}^{*})^{G_{\zero}}$, then the action of $G_{\zero}$ on $\fg_{\one}$ is said to be \emph{polar}. In this case, after fixing a choice of such a $v \in \fg_{\one}$ one can define another detecting subalgebra $\fe = \fe_{\zero} \oplus \fe_{\one}$ via 
	\[
	\fe_{\one}=\left\{ x \in \fg_{\one} : \fg_{\zero}.x \subseteq \fg_{\zero}.v \right\} \quad \text{ and } \quad \fe_{\zero} = [\fe_{\one}, \fe_{\one }].
	\]
Set
	\[
	W = \Norm_{G_{\zero}}(\fe_{\one})/\Stab_{G_{\zero}}(\fe_{\one}).
	\]
It is known for polar actions that $W$ will always be a finite pseduo-reflection group.  For short, we say that $\fg$ is \emph{polar} if the action of $G_{\zero}$ on $\fg_{\one}$ is polar.  Note that both $\ff$ and $\fe$ will again be classical Lie superalgebras and, in the case when the action is both polar and stable, one can choose $x_{0} = v$ and then $\fe \subseteq \ff$, and we will always do so.

By a case-by-case check of the classification of simple classical Lie superalgebras \cite{Kac:1977}, one sees each is either stable or polar, and most are both \cite[Table 5]{Boe:2010}.  The Lie superalgebra $\glmn$ is both stable and polar.  In this case, the detecting subalgebras can be chosen so that $\ff$ is all matrices of the form given in \eqref{E:glmndef} where $A$, $B$, $C$, and $D$ are diagonal with the additional requirement that the $(i,i)$ and $(m+i,m+i)$ entries are equal for all $1 \leq i \leq m$. The subalgebra $\fe \subseteq \ff$ is those matrices with the further condition that $(i,m+i)$ and $(m+i,i)$ are equal for all $1 \leq i \leq m$.

Both $\ff$ and $\fe$ are significantly smaller than $\fg$.  For example, $\fe$ is nearly abelian (namely,  $[\fe_{\zero}, \fe]=0$) and $\ff$ is only a little more complicated.  In many respects, a detecting subalgebra simultaneously plays the role of the Cartan subalgebra in Lie theory and the role of the elementary abelian subgroups of a finite group in support variety theory.  The following theorem illustrates this philosophy.  Note that while the cited result is for classical Lie superalgebras that are stable \emph{and} polar, the proof handles each case separately and so remains valid for those that are stable \emph{or} polar.

\begin{theorem}[{\cite[Theorem 4.1]{Boe:2010}}] \label{T:restriction}
Assume $\fg$ is stable, or polar, or both. Then the restriction maps induced by the inclusions $\ff \hookrightarrow \fg $ and $\fe  \hookrightarrow \fg$
	\begin{align*}
	\res &: \Hbul (\fg, \fg_{\zero}; k) \to \Hbul (\ff, \ff_{\zero}; k), \quad \text{ and } \\
	\res &: \Hbul (\fg, \fg_{\zero}; k) \to \Hbul (\fe, \fe_{\zero}; k)
	\end{align*}
define algebra isomorphisms
	\begin{align*}
	\Hbul (\fg, \fg_{\zero}; k) &\cong \Hbul (\ff, \ff_{\zero}; k)^{N} = S^{\bullet}(\ff^{*}_{\one})^{N}, \quad \text{and} \\
	\Hbul (\fg, \fg_{\zero}; k) &\cong \Hbul (\fe, \fe_{\zero}; k)^{W} = S^{\bullet}(\fe^{*}_{\one})^{W}.
	\end{align*}
\end{theorem}

If $\fg$ is polar, then the previous theorem describes $\Hbul (\fg, \fg_{\zero}; k)$ as the invariant ring of a finite pseduo\-reflection group, and hence $\Hbul (\fg, \fg_{\zero}; k)$ is a polynomial ring.  A case-by-case check verifies that this remains true for the non-polar simple classical Lie superalgebras; see \cite[Appendix]{Boe:2010}.

Given a finite-dimensional Lie superalgebra $\fa$ and a finite-dimensional super\-module $M \in \FF (\fa, \fa_{\zero})$, we write $\VV_{(\fa , \fa_{\zero})}(M)$ for the \emph{relative support variety} defined via the action of the relative co\-hom\-ology ring $\Hbul(\fa,\fa_{\zero};k)$ on the relative cohomology group $\Ext_{\FF(\fa,\fa_{\zero})}^\bullet(M,M) \cong \Hbul(\fa,\fa_{\zero}; M^* \otimes M)$, i.e.,
	\[
	\VV_{(\fa , \fa_{\zero})}(M) = \MaxSpec\big( \Hbul(\fa,\fa_{\zero};k) / \ann_{\Hbul(\fa,\fa_{\zero};k)} \Hbul(\fa,\fa_{\zero}; M^* \otimes M) \big).
	\]
Importantly for both computations and for theoretical results, the elementary structure of $\fe$ allows one to prove that its support varieties have a rank variety description:

\begin{theorem}[{\cite[Theorem 6.4]{Boe:2010}}]
Given an $\fe$-supermodule $M$, set
	\begin{equation}\label{E:RankVariety}
	\VV_{\fe}^{\rank}(M) = \left\{x \in \fe_{\one } :  M \text{ is not projective as an $\subgrp{x}$-supermodule} \right\} \cup \set{0}.
	\end{equation}
Then there is an isomorphism
	\[
	\VV_{\fe}^{\rank}(\C) \cong \VV_{(\fe,\fe_{\zero})}(\C),
	\]
which restricts for each finite-dimensional $\fe$-supermodule $M$ to an isomorphism
	\[
	\VV_{\fe}^{\rank}(M) \cong \VV_{(\fe,\fe_{\zero})}(M)
	\]
\end{theorem}

Let $M$ be a supermodule in $\FF (\fg , \fg_{\zero})$.  If the action of $G_{\zero}$ on $\fg_{\one}$ is stable and/or polar, then there are morphisms of varieties induced by restriction,
	\begin{align*}
	\VV_{(\ff, \ff_{\zero})}(M) &\to \VV_{(\fg, \fg_{\zero})}(M), \\
	\VV_{(\fe, \fe_{\zero})}(M) &\to \VV_{(\fg, \fg_{\zero})}(M).
	\end{align*}
It was conjectured in \cite{Boe:2010} that these maps induce isomorphisms of varieties
	\begin{align}\label{E:SupportIdentifications}
	\VV_{(\ff, \ff_{\zero})}(M)/N &\cong \VV_{(\fg, \fg_{\zero})}(M), \\
	\VV_{(\fe, \fe_{\zero})}(M)/W &\cong \VV_{(\fg, \fg_{\zero})}(M).
	\end{align}
The conjecture was proven for the Type I classical Lie superalgebras (i.e., the classical Lie superalgebras with a compatible $\Z$-grading concentrated in degrees $-1$, $0$, and $1$) in \cite{Lehrer:2011}, and for the remaining simple classical Lie superalgebras in \cite{Grantcharov:2021} using BBW parabolic subalgebras.

One important consequence of this identification is that the rank variety description for  $\fe$-support varieties can be transported to $\ff$- and $\fg$-support varieties.
In particular, this identification implies that support varieties for these Lie superalgebras satisfy the Tensor Product Property:
	\[
	\VV_{(\fg , \fg_{\zero})}(M \otimes N) = \VV_{(\fg , \fg_{\zero})}(M)  \cap \VV_{(\fg , \fg_{\zero})}(N).  
	\]
Support varieties for simple supermodules and other interesting supermodules have been computed in some cases; see \cite{Boe:2009}.

\begin{problem} \label{P:more-examples}
Compute the support varieties of supermodules in other interesting cases.
\end{problem}

\begin{question} \label{Q:realization}
A related question is which subvarieties can be realized as the support variety of a $\fg$-supermodule. As an ingredient in describing the Balmer spectrum of the stable category for $\FF$ (see Section~\ref{SS:BalmerSpectrum}), in \cite{Boe:2017} it was shown for $\glmn$ that every $N$-stable closed subset of $\VV_{(\ff , \ff_{\0})}(\C )$ can be realized as the support variety of a $\fg$-supermodule.  Is this realization statement true for other stable Lie superalgebras?
\end{question}

\subsection{Defect, atypicality, and complexity}\label{SS:defectatypicalityandcomplexity}

Based on other settings, it is reasonable to hope that the dimensions of relative support varieties should encode interesting representation-theoretic invariants.  In many classical settings, the dimension of a module's support variety is equal to the module's \emph{complexity}; that is, the rate of growth of the module's minimal projective resolution.  As we next explain, for classical Lie superalgebras something different occurs.

By definition, a Cartan subalgebra $\fh = \fh_{\zero} \subseteq \fg_{\zero }$ of a classical Lie superalgebra is a choice of a Cartan subalgebra for $\fg_{\zero}$, and the roots $\Phi$ of $\fg$ are the roots of the adjoint action of $\fh$ on $\fg$.  In particular, we write $\Phi_{\zero}$ and $\Phi_{\one}$ for the roots that have a root vector in $\fg_{\zero}$ and $\fg_{\one}$, respectively.  By picking a Borel subalgebra $\fb \subseteq \fg$, there is a corresponding set of postive (resp., negative) roots $\Phi^{+}$ (resp., $\Phi^{-}$).  For $r \in \Z_{2}$, set $\Phi^{\pm}_{r}= \Phi^{\pm} \cap \Phi_{r}$. Put 
	\[
	\rho = \frac{1}{2}\sum_{\alpha \in \Phi^{+}_{\zero}} \alpha - \frac{1}{2}\sum_{\alpha \in \Phi^{+}_{\one}} \alpha.
	\]

A classical Lie superalgebra is called \emph{basic} if it admits a nondegenerate, invariant, super\-sym\-metric, even bilinear form $(-,-)$.  For example, the supertrace defines such a bilinear form on $\glmn$ by the formula $(X,Y) = \strace (XY)$.  By the usual identifications, a basic classical Lie super\-algebra then has a nondegenerate bilinear form on $\fh^{*}$, which we also write as $(-,-)$.  By definition, the \emph{defect} of $\fg$, denoted $\defect (\fg )$, is the cardinality of a maximal set of mutually orthogonal, isotropic elements of $\Phi_{\one}^{+}$ with respect to this bilinear form on $\fh^*$.  For example, the defect of $\glmn$ is the minimum of $m$ and $n$.  We will assume $\fg$ is basic classical for the remainder of this section.

Recall that for a simple basic classical Lie superalgebra the cohomology ring $\Hbul(\fg , \fg_{\0}; k)$ is a polynomial ring in some number of variables.  A case-by-case check verifies that the number of variables is precisely the defect of $\fg$:
	\[
	\defect (\fg ) = \dim \Hbul (\fg , \fg_{\zero};k) = \dim \VV_{(\fg , \fg_{\zero})}(k).
	\]
Thus the dimension of the ambient space for support varieties encodes the defect of $\fg$. This raises the question of whether root system combinatorics can also be used to compute the dimension of support varieties for interesting supermodules (e.g., the simple supermodules). As we next explain, the answer appears to be yes.

Given $\lambda \in \fh^{*}$, the \emph{atypicality} of $\lambda$, denoted $\atyp (\lambda)$, is the cardinality of a maximal set of mutually orthogonal, isotropic elements of $\Phi_{\one}^{+}$ that are orthogonal to $\lambda + \rho$. Clearly, 
	\[
	\atyp (\lambda) \leq \defect (\fg ).
	\]

As for ordinary Lie algebras, the simple $\g$-supermodules in $\FF$ can be labeled by their highest weight with respect to our choice of Cartan and Borel subalgebras.  Let 
	\[
	\left\{L(\lambda) : \lambda \in X^{+} \subseteq \fh^{*} \right\}
	\]
be a complete, irredundant set of simple $\g$-supermodules in $\FF$. It is well-known that the atypicality of the highest weight encodes important representation-theoretic information about the simple supermodule it labels.  For example, Kac and Wakimoto \cite{KacWakimoto} conjectured that the superdimension of a simple supermodule, 
	\[
	\sdim_{k} L(\lambda) := \dim_{k}L(\lambda)_{\0}- \dim_{k}L(\lambda)_{\1},
	\]
should be nonzero if and only if $\atyp (\lambda) = \defect (\fg )$.  In \cite{KGPM}, the second author along with Geer and Paturau-Mirand formulated a generalization of this conjecture for all degrees of atypicality by relating it to the nonvanishing of certain modified dimension functions.  The Kac--Wakimoto conjecture and its generalization are now known to be true for $\glmn$ and $\osp (m|2n)$, i.e., for types $ABCD$ in the Kac classification; see \cite{SerganovaGenKW,KujawaGenKW}.

It is easy to see from the rank variety description given by \eqref{E:RankVariety} and \eqref{E:SupportIdentifications} that if $\sdim_{k} L(\lambda) \neq 0$, then $\VV_{(\fg, \fg_{\0}) }(L(\lambda)) = \VV_{(\fg , \fg_{\0})}(k)$.  That is, the nonvanishing of the superdimension is related to the support variety being as large as possible, and hence having dimension equal to the defect of $\fg$.  In light of the generalized Kac--Wakimoto conjecture, one could imagine that atypicality is related to the dimension of the support variety more generally.  The next result confirms this speculation.

\begin{theorem}[{\cite[Theorem 4.7]{Boe:2009}, \cite[Theorem 4.3]{KujawaGenKW}}]
Let $\fg$ be $\glmn$ or $\osp (m|2n)$, and let $L(\lambda)$ be a simple $\g$-supermodule in $\FF$. Then,
	\[
	\dim \VV_{(\fg, \fg_{\0})}(L(\lambda))=\atyp (\lambda).
	\]
\end{theorem}

This gives a cohomological/geometric interpretation of atypicality. Moreover, it naturally suggests to define the atypicality of an arbitrary supermodule in $\FF$ as the dimension of its support variety.

\begin{question}
Suppose $\fg$ is a classical Lie superalgebra that is not basic, e.g., a Lie superalgebra of type $P$ or type $Q$. Can one give a combinatorial method of computing the dimension of a simple $\g$-supermodule's support variety?
\end{question}

The dimension of a relative support variety does not directly equal the complexity of a supermodule, but we can ask if they are nevertheless related.  It turns out that there are two relevant notions of complexity.  As mentioned earlier, the \emph{complexity} of a supermodule $M \in \FF$ is the rate of growth of the dimensions of the terms in a minimal projective resolution $P_\bullet$ in $\FF$ of $M$. That is, the complexity of $M$ is the minimal integer $c$ such that there is a constant $K > 0$ for which $\dim_{k} P_{d} \leq Kd^{c-1}$ for all $d \geq 0$ (declaring the complexity to be infinite if no such $c$ exists).  A second invariant is the \emph{$z$-complexity} of $M$, defined to be the rate of growth of the number of indecomposable summands of the minimal resolution.  Importantly, unlike complexity, $z$-complexity is invariant under category equivalences.  Let us write $c_{\FF}(M)$ and $z_{\FF}(M)$ for the complexity and $z$-complexity of $M$, respectively.

The \emph{odd nullcone} of a Lie superalgebra $\g$ is the set
	\[
	\calN_{\odd}(\fg ) = \Chi_{\fg}(k) = \set{ x \in \fg_{\one} : [x,x]=0}.
	\]
Given $M$ in $\FF$, one has the \emph{associated variety} introduced by Duflo and Serganova \cite{Duflo:2005}:
	\begin{equation} \label{eq:DF-variety-char0}
	\Chi_{\fg }(M) = \set{ x \in \calN_{\odd}(L) : M|_{\subgrp{x}} \text{ is not free}} \cup \set{0}.
	\end{equation}
If $\fg$ is $\glmn$, $\osp (2|2n)$, $\osp (k|2)$, $D(2,1;\alpha)$, $G(3)$, or $F(4)$, then it is known that the complexity and $z$-complexity of a simple supermodule each have a geometric interpretation.

\begin{theorem}[\cite{Boe:2012,El-Turkey:2016,El-Turkey:2018}]\label{T:complexity}
Let $\fg$ be a Lie superalgebra from the above list and let $L$ be a simple supermodule in $\FF$. Then 
	\begin{align*}
	c_{\FF }(L) &= \dim \VV_{(\fg , \fg_{\zero})}(L) + \dim \Chi_{\fg}(L), \quad \text{and} \\
	z_{\FF}(L) &= \dim \VV_{(\ff , \ff_{\zero})} (L).
	\end{align*}
\end{theorem}

If $\fg$ is either $\glmn$ or $\osp (2|2n)$, then $\fg$ is of Type I and there is a compatible $\Z$-grading $\fg = \fg_{-1} \oplus \fg_{0} \oplus \fg_{1}$.  Given a simple $\fg_{0}$-supermodule $S$, one makes it into a $\fg_{0} \oplus \fg_{1}$-supermodule by having $\fg_{1}$ act trivially, and one then defines the associated \emph{Kac supermodule}
	\[
	K(S) = U(\fg ) \otimes_{U(\fg_{0} \oplus \fg_{1})} S.
	\]
In the cases $\g = \glmn$ and $\g = \osp (2|2n)$, the previous theorem also holds true if $L$ is replaced with $K(S)$.  Similarly, the type $P$ Lie superalgebra has so-called \emph{thin} Kac supermodules, and the above formulas hold for these as well \cite{Boe:2020}.

Unfortunately, the above results are obtained by determining explicit formulas for the complexity and $z$-complexity and, separately, the dimensions of the relevant varieties. The proof is then completed by nodding sagely and observing that the two numbers in each equality are indeed equal. This leads to the following important open question:

\begin{question}
Do the formulas for $c_{\FF}(L)$ and $z_{\FF}(L)$ given in Theorem~\ref{T:complexity} hold for arbitrary classical Lie superalgebras and arbitrary supermodules in $\FF$?  
\end{question}

Answering this question will no doubt require new insights.  In particular, the Duflo--Serganova associated variety is not (yet) known to have a cohomological description in characteristic zero. The fact that it appears alongside the support variety in these complexity formulas remains unexplained.  For Type I simple Lie superalgebras, a hint that it has a role to play can be found in the fact that the associated variety vanishes if and only if the supermodule in question is projective \cite{Gorelik:2022}.

It is worth remarking that an interesting theory of relative support varieties can also be developed for certain non-classical Lie superalgebras.  If $\fg = \bigoplus_{k \in \Z} \fg_{k}$ has a suitable $\Z$-grading, then one can instead consider relative cohomology for the pair $(\fg , \fg_{0})$.  This setup applies to the simple Lie superalgebras of type $WSH$ in the Kac classification. For example, in \cite{Bagci:2008} the theory is developed for the Lie superalgebras $W(n)$ and, once again, the support variety of a simple supermodule is shown to capture its atypicality.

\subsection{The Balmer Spectrum} \label{SS:BalmerSpectrum}

Let $\KK := \Stab (\FF )$ be the stable supermodule category for $\FF$ obtained by factoring out the homomorphisms that factor through a projective supermodule.  Since the projective and injective supermodules in $\FF$ coincide, $\KK$ is a triangulated category.  The tensor structure inherited from $\FF$ makes $\KK$ a symmetric tensor triangulated category.  To such a category, Balmer \cite{Balmer:2005} introduced the notion of the spectrum of $\KK$, denoted $\Spec (\KK )$.

Briefly, a \emph{(thick) tensor ideal} of $\KK$ is a full, replete, triangulated subcategory $\mathcal{I}$ of $\KK$ that satisfies the following conditions:
	\begin{itemize}
	\item If $X$ and $Y$ are objects of $\KK$ and $X \oplus Y$ is an object of $\mathcal{I}$, then $X$ and $Y$ are objects of $\mathcal{I}$;
	\item If $X$ is an object of $\KK$ and $Y$ is an object of $\mathcal{I}$, then $X \otimes Y$ is an object of $\mathcal{I}$.
	\end{itemize}
A tensor ideal $\mathcal{I}$ is \emph{prime} if it further satisfies the condition that whenever $X \otimes Y$ is an object of $\mathcal{I}$, then either $X$ or $Y$ is an object of $\mathcal{I}$.

The Balmer spectrum $\Spec (\KK )$ is the collection of all proper prime tensor ideals in $\KK$, considered as a topological space via the Zariski topology. The \emph{support} of an object $X$ in $\KK$ is defined by
	\[
	\supp (X) = \left\{\mathcal{P} \in \Spec (\KK ) : X \text{ is not an object of } \mathcal{P} \right\}.
	\]
The Balmer spectrum and support are known to have many of the desirable properties of a support variety theory. They satisfy a certain universal property for such support variety theories, and they can be used to classify the tensor ideals of $\KK$.  For classical Lie superalgebras we have seen four interrelated support theories: the support varieties of $\fg$, $\ff$, $\fe$, and the associated variety of Duflo and Serganova.  It is an obvious question to ask how these are related to Balmer's `universal' support theory, and to hopefully give a more concrete description of Balmer's theory for $\Stab (\FF )$.

As we explain, the detecting subalgebra $\ff$  plays a key role.  First, recall that $N$ is the normalizer of the action of $G_{\zero}$ on $\ff_{\1}$.  There is an action of $N$ on $S^{\bullet}(\ff_{\1}) = \Hbul (\ff , \ff_{\zero};  k)$.  Let $X = \Proj( S^\bullet(\ff_{\one}))$ be the set of homogeneous (with respect to the cohomological $\Z$-grading) prime ideals of $S^{\bullet}(\ff_{\1})$, and let $X_{N} = \NProj(S^\bullet(\ff_{\one})) \subseteq X$ be the set of homogeneous (with respect to the $\Z$-grading) $N$-prime ideals of $S^{\bullet}(\ff_{\1})$.  By an $N$-prime ideal we mean an $N$-stable ideal of $S^{\bullet}(\ff_{\1}^{*})$ that is prime among the set of all $N$-stable ideals; see for example \cite{Lorenz:2009}. There is a canonical map $\pi: X \to X_{N}$ given by $P \mapsto \bigcap_{n \in N} n \cdot P$.

The key feature of $X_{N}$ is that its closed sets are in bijection, via $\pi$, with the $N$-stable closed sets of $X$. In particular, given a $\fg$-supermodule $M$ in $\FF$, the projectivization of the support variety $\VV_{(\ff, \ff_{\zero})}(M)$ is an $N$-stable closed set in $X$, and its image under $\pi$ is a closed set in $X_{N}$.  Recall that in a topological space a subset is called \emph{specialization closed} if it is a union of (arbitrarily many) closed subsets.  Then the main result of \cite{Boe:2017} is the following:

\begin{theorem}[{\cite[Theorems 5.2.1~and~5.2.2]{Boe:2017}}]
Let $\fg$ be $\glmn$ or $\osp (2|2n)$.  Then there is a homeomorphism 
	\[
	X_{N} \to \Stab (\FF (\fg , \fg_{\0}))
	\]
that identifies the cohomological support variety (with respect to the detecting subalgebra $\ff$) with the Balmer support.

Moreover, the thick tensor ideals of $\Stab (\FF (\fg , \fg_{\zero}) )$ are in bijection with the specialization closed subsets of $X_{N}$.
\end{theorem}

\begin{question}
Can one obtain a similar description of the Balmer spectrum for classical Lie superalgebras in other types? 
\end{question}

If one has a support variety theory, then there are results that can be used to identify it with Balmer's support theory and, in turn, to identify the Balmer spectrum; see \cite[Theorem 3.5.1]{Boe:2017} or \cite[Theorem 1.5]{DellAmbrogio:2010}. Using these results requires that the support variety theory in question satisfy a list of properties, including the `realizability' condition that whenever $W$ is a closed subset of the ambient space, then there is an object $M$ in $\KK$ for which the support theory's value on $M$ is $W$.  In the context of the previous theorem, this says, roughly, that given a closed $N$-stable subset of $\Proj (S^{\bullet}(\ff_{\1}))$, we must find a finite-dimensional $\fg$-supermodule that has the given closed subset as its $\ff$-support variety.  The standard realization arguments are not helpful here.  Instead, the necessary closed sets are constructed `by hand' using Kac supermodules.  Since Kac supermodules only exist for Type I Lie superalgebras, new methods would be required for the other types.  Thus, Question~\ref{Q:realization}.

\section{Support varieties for finite supergroup schemes in characteristic zero} \label{sec:supergroup-char-0}

\subsection{Finite supergroup schemes in characteristic zero}

The data of a finite supergroup scheme is equivalent to the data of a finite-dimensional cocommutative Hopf superalgebra. Over an algebraically closed field of characteristic zero, these have a particularly simple form:

\begin{theorem}[{\cite[Corollary 3.1.2]{Drupieski:2019a}}]
Let $k$ be an algebraically closed field of characteristic zero, and let $A$ be a finite-dimensional cocommutative Hopf superalgebra over $k$. Then there exists a finite group $G$, a finite-dimensional odd superspace $V$, and a representation of $G$ on $V$ such that $A$ is isomorphic as a Hopf superalgebra to the smash product algebra $\Lambda(V) \# kG$. Here $kG$ denotes the group algebra of $G$ over $k$, considered as a purely even superalgebra.
\end{theorem}

We denote the finite supergroup scheme corresponding to the algebra $\Lambda(V) \# kG$ by $V \rtimes G$, and we write $V/G$ for the set of $G$-orbits in $V$. Over a field of characteristic zero, the group algebra $kG$ is semisimple. This leads to the following calculation of the cohomological spectrum of $V \rtimes G$:

\begin{theorem}[{\cite[Theorem 3.2.2]{Drupieski:2019a}}]
Let $G$ be a finite group, and let $V$ be a finite-dimensional purely odd $kG$-module. Then there exist isomorphisms of varieties
	\[
	\abs{V \rtimes G} \cong \MaxSpec( S(V^*)^G) \cong V/G.
	\]
\end{theorem}

Write $[v]$ for the $G$-orbit of an element $v \in V$, and write $\subgrp{v}$ for the $k$-subalgebra of $\Lambda(V)$ generated by $v$. By convention, $\subgrp{0} = k$. The category of $V \rtimes G$-supermodules is equivalent to the category of supermodules over the superalgebra $\Lambda(V) \# kG$. Now one gets the following rank variety description for the support varieties of finite-dimensional $V \rtimes G$-supermodules:

\begin{theorem}[{\cite[Theorem 3.2.3]{Drupieski:2019a}}]
Let $G$ be a finite group, and let $V$ be a finite-dimensional purely odd $kG$-module. Let $M$ be a finite-dimensional $V \rtimes G$-supermodule. Then
	\[
	\abs{V \rtimes G}_M \cong \set{ [v] \in V/G: M|_{\subgrp{v}} \text{ is not free}} \cup \set{0}.
	\]
\end{theorem}

\hypersetup{bookmarksdepth=-1}
\part{Positive characteristic} \label{part:char-p}
\hypersetup{bookmarksdepth=3}

In Part \ref{part:char-p}, let $k$ be a field of characteristic $p > 0$. Usually $p \neq 2$, although some results can be extended to the case $p=2$; see \cite{Drupieski:2022b}. For convenience we will assume that $k$ is algebraically closed, although some results may hold without this assumption. Given a $k$-vector space $V$ and an integer $r \geq 1$, let $V^{(r)} = V \otimes_{\varphi_r} k$ be the $r$-th Frobenius twist of $V$, i.e., the $k$-vector space obtained from $V$ via base change along the ($r$-th) Frobenius endomorphism $\varphi_r: k \to k$, $\lambda \mapsto \lambda^{p^r}$. If $W = V^{(r)}$, then we may write $V = W^{(-r)}$. More generally, if $X$ is an affine $k$-scheme (resp.\ affine algebraic variety) with coordinate algebra $k[X]$, then we may write $X^{(r)}$ for the scheme (resp.\ variety) with coordinate algebra $k[X^{(r)}] = k[X]^{(r)}$, and if $Y = X^{(r)}$, then we may write $X = Y^{(-r)}$.

\section{Support varieties via ordinary Lie super\-algebra cohomology} \label{sec:support-ordinary}

\subsection{Lie super\-algebra cohomology in positive characteristic}

Let $L$ be a finite-dimensional Lie super\-algebra over the field $k$. As discussed in Section \ref{subsec:generalities-cohomology}, the Lie super\-algebra cohomology ring $\Hbul(L,k)$ can be computed as the cohomology of the superexterior algebra $\bsl(L^*)$ with respect to the Koszul differential $\partial$. Since $\partial$ makes $\bsl(L^*)$ into a differential graded super\-algebra, and since $\bsl(L^*)$ is graded-commutative in the sense of \eqref{eq:graded-commutative}, it follows for all $f \in \Lone^* \subseteq \bsl^1(L^*)$ that $\partial(f^p) = 0$. Then there exists a map of graded super\-algebras
	\begin{equation} \label{eq:varphi}
	\varphi: S(\Lone^*[p])^{(1)} \to \Hbul(L,k)
	\end{equation}
induced by the $p$-power map on $S(\Lone^*)$ and the inclusion $S(\Lone^*) \subset \bsl(L^*)$. If $M$ and $N$ are finite-dimensional $L$-supermodules, then $\Ext_L^\bullet(M,N)$ is a finite module over the image of $\varphi$. In particular, $\Hbul(L,k)$ is a finitely-generated $k$-algebra, and $\Ext_L^\bullet(M,N)$ is a finite $\Hbul(L,k)$-module.

Given an $L$-supermodule $M$, let $J_L(M)$ be the kernel of the composite map
	\[
	S(\Lone^*)^{(1)} \xrightarrow{\varphi} \Hbul(L,k) \xrightarrow{- \otimes M} \Ext_L^\bullet(M,M),
	\]
and set
	\[\Chi_L(M) = \MaxSpec\big( S(\Lone^*)^{(1)}/J_L(M) \big).
	\]
Then $\varphi$ induces for each $L$-supermodule $M$ a homeo\-morphism of varieties
	\begin{equation} \label{eq:U(L)-ChiLM}
	\varphi^*: \abs{U(L)}_M \simeq \Chi_L(M).
	\end{equation}
In the case $M = k$, this provides the following explicit description of the cohomological spectrum for $U(L)$ in terms of the \emph{odd nullcone} of $L$.

\begin{theorem}[{\cite[Theorem 4.2.4]{Drupieski:2019a}}]
Let $L$ be a finite-dimensional Lie super\-algebra over $k$. Then
	\[
	\Chi_L(k)^{(-1)} = \calN_{\odd}(L) := \set{ x \in \Lone : [x,x] = 0}.
	\]
\end{theorem}

\subsection{Support varieties as rank varieties}

One of the main results of \cite{Drupieski:2021b} was the following `rank variety' description for $\Chi_L(M)$ when $M$ is finite-dimensional:

\begin{theorem}[{\cite[Corollary 3.3.2]{Drupieski:2021b}}]
Let $L$ be a finite-dimensional Lie super\-algebra over $k$, and let $M$ be a finite-dimensional $L$-supermodule. Then
	\begin{equation} \label{eq:L-rank-equals-support}
	\Chi_L(M)^{(-1)} = \set{ x \in \calN_{\odd}(L) : M|_{\subgrp{x}} \text{ is not free}} \cup \set{0}.
	\end{equation}
\end{theorem}

Here $M|_{\subgrp{x}}$ denotes the restriction of $M$ to the $k$-subalgebra of $U(L)$ generated by $x$; for $x = 0$ this subalgebra is just the field $k$ (and hence $M|_{\subgrp{x}}$ is free), while for $0 \neq x \in \calN_{\odd}(L)$ this subalgebra is of the form $k[x]/(x^2)$. The rank variety on the right-hand side of \eqref{eq:L-rank-equals-support} is equivalent in definition to the `associated variety' of an $L$-supermodule originally defined---without reference to cohomology---by Duflo and Serganova for Lie super\-algebras in characteristic zero \cite{Duflo:2005}, and mentioned earlier in Section~\ref{SS:defectatypicalityandcomplexity}.

Let us say a few words about the strategy used in \cite{Drupieski:2021b} to prove \eqref{eq:L-rank-equals-support}. Let $\Chi_L'(M)$ be the closed subset of $\Chi_L(k) = \calN_{\odd}(L)$ such that $\Chi_L'(M)^{(-1)}$ is the right-hand side of \eqref{eq:L-rank-equals-support}. First, using naturality it is relatively easy to show that $\Chi_L'(M) \subseteq \Chi_L(M)$. Equality for arbitrary $L$ and $M$ can then be reduced to showing for all $m \geq 1$ that equality holds for the general linear Lie super\-algebra $\gl(m|m)$ and its natural representation $M = k^{m|m}$. To prove equality for $\g = \gl(m|m)$, we consider the Clifford filtration on $\g$, that is, the Lie super\-algebra filtration $F^0 \g \subseteq F^1 \g \subseteq F^2 \g$ defined by $F^0 \g = 0$, $F^1 \g = \gone$, and $F^2 \g = \g$. The associated graded Lie super\-algebra $\wtg$ is concentrated in $\Z$-degrees $1$ and $2$, with $\wtg_1 \cong \gone$ and $\wtg_2 \cong \gzero$ as superspaces. The Lie bracket on $\wtg_1$ identifies with the original Lie bracket on $\gone$ (and hence retains information about the odd nullcone of $\g$), while $\wtg_2$ is now central in $\wtg$. Given a choice of generators for the $\g$-supermodule $M$, one can also define an associated graded $\wtg$-supermodule $\wt{M}$. Now the Clifford filtration gives rise to a spectral sequence relating $\Hbul(\wtg,k)$ and $\Hbul(\g,k)$, and by studying this spectral sequence one can relate the support of $M$ as a $\g$-supermodule and the support of $\wt{M}$ as a $\wtg$-supermodule. Finally, the support of $\wt{M}$ as a $\wtg$-supermodule is calculable, and thus can be used to help put an upper bound on the size of $\Chi_\g(M)$. The support of the $\wtg$-supermodule $\wt{M}$ is easier to calculate partly because $\wtg$ has a simpler Lie super\-algebra structure than $\wtg$, and partly because the $\wtg$-action on $\wt{M}$ factors through the action of a $p$-nilpotent restricted Lie superalgebra (equivalently, through the action of an infinitesimal unipotent supergroup scheme), for which we have a clearer picture of the support theory; see Remark \ref{rem:p-nilpotent-explicit}.

An immediate consequence of the `rank variety' description \eqref{eq:L-rank-equals-support} and the representation theory of the Hopf algebra $k[x]/(x^2)$ is that support varieties satisfy the Tensor Product Property:

\begin{proposition} \label{P:TensorProductProperty}
Let $L$ be a finite-dimensional Lie super\-algebra over a field $k$ of characteristic $p \geq 3$, and let $M$ and $N$ be finite-dimensional $L$-supermodules. Then
	\[
	\Chi_L(M \otimes N) = \Chi_L(M) \cap \Chi_L(N).
	\]
Consequently, $\abs{U(L)}_{M \otimes N} = \abs{U(L)}_M \cap \abs{U(L)}_N$.
\end{proposition}

As mentioned in Section \ref{subsec:basic-notions-support}, if $A$ is a \emph{finite-dimensional} (and hence self-injective) Hopf super\-algebra whose cohomology satisfies the finiteness properties \ref{FG-ring} and \ref{FG-modules}, then for each finite-dimensional $A$-super\-module $M$ one gets $\dim(\abs{A}_M) = \cx_A(M)$. In particular, one gets that $\abs{A}_M = \set{0}$ if and only if $M$ is projective. The cohomology of $U(L)$ does satisfy \ref{FG-ring} and \ref{FG-modules}, but $U(L)$ is not finite-dimensional unless $L$ is a purely odd (abelian) Lie super\-algebra. While $U(L)$ is not self-injective in general, it is Gorenstein (i.e., it has finite injective dimension as a module over itself), and a finitely-generated $U(L)$-supermodule $M$ has finite projective dimension if and only if it has finite injective dimension; see \cite[\S2.4]{Drupieski:2021b}. For a finite-dimensional $U(L)$-supermodule $M$, one can show that $\abs{U(L)}_{M} = \set{0}$ if and only if $M$ is of finite projective dimension. Reformulated in terms of rank varieties, this yields the following statement:

\begin{theorem}[{\cite[Theorem 3.4.2, Corollary 3.4.3]{Drupieski:2021b}}] \label{thm:fin-proj-dim}
Let $L$ be a finite-dimensional Lie super\-algebra over an algebraically closed field $k$ of characteristic $p \geq 3$, and let $M$ be a finite-dimensional $L$-supermodule. Then
	\[
	\set{ x \in \calN_{\odd}(L) : M|_{\subgrp{x}} \text{ is not free}} = \emptyset \quad \text{if and only if} \quad \projdim_{U(L)}(M) < \infty.
	\]
In particular, if $\calN_{\odd}(L) = \set{0}$, then $U(L)$ has finite global dimension (and conversely).
\end{theorem}

For the converse of the last statement, one observes that $\gldim(U(L')) \leq \gldim(U(L))$ for any Lie sub-super\-algebra $L'$ of $L$, as a consequence of the fact that $U(L)$ is free over $U(L')$. So if $U(L)$ has finite global dimension, it cannot have any subalgebras of the form $k[x]/(x^2)$ for $x \in \Lone$. The last statement of Theorem \ref{thm:fin-proj-dim} was previously known in characteristic zero by work of B{\o}gvad \cite{Bo-gvad:1984}; for details, see \cite[Theorem 17.1.2]{Musson:2012}.

\begin{problem}
Let $L$ be a finite-dimensional Lie super\-algebra over $k$, and let $M$ be a finite-dimensional $L$-supermodule. Provide a general representation-theoretic interpretation for the geometric dimension of the support variety $\abs{U(L)}_M \simeq \Chi_L(M)$.
\end{problem}

\section{Finite group schemes: recollections from the non-super theory} \label{sec:non-super}

\subsection{Restricted Lie algebras}

Let $\g$ be a finite-dimensional restricted Lie algebra over $k$. The restricted enveloping algebra of $\g$, denoted $V(\g)$, is a finite-dimensional cocommutative Hopf algebra over $k$. The dual Hopf algebra $V(\g)^* = \Hom_k(V(\g),k)$ is then a finite-dimensional commutative Hopf algebra over $k$ with the property that $f^p = 0$ for each element $f$ of the augmentation ideal $I_\ve$ of $V(\g)^*$. Thus $V(\g)^*$ is the coordinate algebra $k[G]$ of a height-$1$ infinitesimal group scheme $G$, i.e., a finite $k$-group scheme that is equal to its own first Frobenius kernel $G_{(1)}$. The category of (left) $V(\g)$-modules is then equivalent to the category of rational (left) $G$-modules, and cohomology for $V(\g)$ may be identified with (rational) cohomology for $G$. For more details, see \cite[I.8, I.9]{Jantzen:2003}.

Let $M$ be a $V(\g)$-module (equivalently, a rational $G$-module). The powers of the augmentation ideal $I_\ve \subset k[G]$ give rise to a filtration on the Hochschild complex $\Cbul(G,M)$ that computes the cohomology group $\Hbul(G,M) = \Hbul(V(\g),M)$. The filtration gives rise in turn to a spectral sequence, which for $p > 2$ can be written in the form
	\begin{equation} \label{eq:FPspecseq}
	E_2^{i,j}(M) = S^{i/2}(\g^*)^{(1)} \otimes \opH^j(\g,M) \Rightarrow \opH^{i+j}(V(\g),M).
	\end{equation}
Here $\Hbul(\g,M)$ is the (ordinary) Lie algebra cohomology group of $M$, and the superscript $i/2$ means that $E_2^{i,j}(M) = 0$ for odd $i$. In the case $M = k$, the horizontal edge map of \eqref{eq:FPspecseq} defines a finite map of graded $k$-algebras
	\begin{equation} \label{eq:Phi*}
	\Phi^*: S(\g^*[2])^{(1)} \to \Hbul(V(\g),k).
	\end{equation}
Here $S(\g^*[2])$ means that we consider the symmetric algebra $S(\g^*)$ as generated in cohomological degree $2$. Passing to maximal ideal spectra, one gets a finite morphism of varieties
	\[
	\Phi: \abs{V(\g)} := \MaxSpec\big( \Hbul(V(\g),k) \big) \to \MaxSpec\big( S(\g^*[2])^{(1)} \big) = \g^{(1)}
	\]

Jantzen \cite{Jantzen:1986} showed, as a consequence of explicit calculations in the case $\g = \gl_n$, that the image of the morphism $\Phi$ is the restricted nullcone of $\g$:
	\begin{equation} \label{eq:image-spectrum-nullcone}
	\Phi( \abs{V(\g)} )^{(-1)} = \calN_1(\g) := \{ x \in \g : x^{[p]} = 0 \}.
	\end{equation}
Then for each finite-dimensional $V(\g)$-module $M$, Friedlander and Parshall \cite{Friedlander:1986b} computed the image of the support variety $\abs{V(\g)}_M$ under the morphism $\Phi$:
	\begin{equation} \label{eq:restricted-support-equals-rank}
	\Phi( \abs{V(\g)}_M )^{(-1)} = \set{ x \in \calN_1(\g) : M|_{\subgrp{x}} \text{ is not free}} \cup \set{0}.
	\end{equation}
Here $M|_{\subgrp{x}}$ denotes the restriction of $M$ to the $k$-subalgebra of $V(\g)$ generated by $x$. If $x = 0$ this subalgebra is just the field $k$, while for $x \neq 0$ the subalgebra has the form $k[x]/(x^p)$.

\subsection{Cohomology of infinitesimal group schemes} \label{sec:cohomology-infinitesimal}

For this section, fix an integer $r \geq 1$.

Using the theory of strict polynomial functors, Friedlander and Suslin \cite{Friedlander:1997} established the existence of certain `universal extension classes' for the rational cohomology of $GL_n$:
	\begin{equation} \label{eq:universal-classes}
	e_i^{(r-i)} \in \opH^{2p^{i-1}}(GL_n,\gl_n^{(r)}), \quad 1 \leq i \leq r.
	\end{equation}
Here $\gl_n$ denotes the adjoint representation of $GL_n$, and the superscript $(r)$ indicates both that the vector space structure of $\gl_n$ has been twisted by the $r$-th Frobenius morphism of the field, and also that the rational $GL_n$-module structure has been twisted via pre-composition with the $r$-th Frobenius morphism $F^r: GL_n \to (GL_n)^{(r)}$ of the scheme $GL_n$. The restriction of $e_i^{(r-i)}$ to the $r$-th Frobenius kernel $\GLnr$ of $GL_n$ determines an element of
	\[
	\opH^{2p^{i-1}}(\GLnr,\gl_n^{(r)}) \cong \opH^{2p^{i-1}}(\GLnr,k) \otimes \gl_n^{(r)} \cong \Hom_k(\gl_n^{*(r)},\opH^{2p^{i-1}}(\GLnr,k)),
	\]
i.e., it determines a linear map $\gl_n^{*(r)} \to \opH^{2p^{i-1}}(\GLnr,k)$, which then extends multiplicatively to a homomorphism of graded $k$-algebras $S(\gl_n^{*(r)}[2p^{i-1}]) \to \Hbul(\GLnr,k)$. Taking the product of these homomorphisms, one gets a map of graded $k$-algebras
	\begin{equation} \label{eq:phi-GLnr} \textstyle
	\phi_{\GLnr} : \bigotimes_{i=1}^r S( \gl_n^{*(r)}[2p^{i-1}]) \to \Hbul(\GLnr,k).
	\end{equation}

Now let $G$ be an infinitesimal $k$-group scheme of height $\leq r$. Thus $G$ is a finite $k$-group scheme with the property that $f^{p^r} = 0$ for each element $f$ of the augmentation ideal of $k[G]$. Equivalently, $G$ is equal to its own $r$-th Frobenius kernel $G_{(r)}$. Set $\g = \Lie(G)$, and fix for some $n$ a choice of closed embedding $\iota: G \hookrightarrow GL_n$, i.e., an embedding of $G$ as a closed subgroup scheme of $GL_n$. Then the image of $G$ is contained in $\GLnr$, and $\iota$ differentiates to an injective homomorphism $d\iota: \g \hookrightarrow \gl_n$ of restricted Lie algebras. Composing $\phi_{\GLnr}$ with the restriction map in cohomology, one gets a homomorphism of graded $k$-algebras
	\[ \textstyle
	\phi_G: \bigotimes_{i=1}^r S( \gl_n^{*(r)}[2p^{i-1}]) \to \Hbul(G,k).
	\]
One of the main results of \cite{Friedlander:1997} then states that $\phi_G$ is a finite algebra map, and that $\Hbul(G,M)$ is finite over the image of $\phi_G$ whenever $M$ is a finite-dimensional rational $G$-module. In particular, $\Hbul(G,k)$ is a finitely-generated $k$-algebra, and $\Hbul(G,M)$ is a finite $\Hbul(G,k)$-module. In the special case $r=1$, the map $\phi_G$ reduces to the map \eqref{eq:Phi*} considered by Friedlander and Parshall.

\subsection{Support varieties for infinitesimal group schemes} \label{sec:SFB-support-varieties}

Given an affine $k$-group scheme $G$ with coordinate algebra $k[G]$, and given a commutative $k$-algebra $A$, let $G_A = G \otimes_k A$ be the affine $A$-group scheme with coordinate algebra $A[G_A] := k[G] \otimes_k A$. If $M$ is a rational $G$-module, then $M_A := M \otimes_k A$ is a rational $G_A$-module.

Let $\Ga$ be the additive group scheme, and $\Gar$ its $r$-th Frobenius kernel. Write $\CAlg(k)$ for the category of commutative (unital, associative) $k$-algebras. The functor
	\[
	V_r(G) : \CAlg(k) \to \Set
	\]
is defined by
	\begin{equation} \label{eq:Vrg-non-super}
	V_r(G)(A) = \Hom_{\Grp/A}(\Gar \otimes_k A, G \otimes_k A),
	\end{equation}
the set of $A$-group scheme homomorphisms $\nu: \Gar \otimes_k A \to G \otimes_k A$. If $G$ is an affine algebraic $k$-group scheme (i.e., if $k[G]$ is a finitely-generated $k$-algebra), then by \cite[Theorem 1.5]{Suslin:1997}, $V_r(G)$ admits the structure of an affine $k$-scheme of finite type, called the scheme of one-parameter subgroups of height $\leq r$ in $G$. In the special case $r=1$, $V_1(G)$ identifies with $\calN_1(\g)$, the restricted nullcone of $\g = \Lie(G)$. More generally, for $G = \GLnr$ one gets a natural identification
	\begin{equation} \label{eq:Vr(GLn)}
	V_r(\GLnr)(A) = \set{ (\alpha_0,\cdots,\alpha_{r-1}) \in \gl_n(A) : \alpha_i^p = 0 = [\alpha_j,\alpha_\ell] \text{ for all } 0 \leq i,j,\ell < r}.
	\end{equation}

The `group algebra' of a finite $k$-group scheme, denoted $kG$, is the finite-dimensional Hopf algebra that is dual to the coordinate algebra of $G$: $kG = k[G]^*$. In the notation of \cite[I.8]{Jantzen:2003}, $kG = M(G)$, and if $G$ is infinitesimal, then $kG = \Dist(G)$. For a finite $k$-group scheme, the category of rational (left) $G$-modules is equivalent to the category of (left) $kG$-modules \cite[I.8.6]{Jantzen:2003}, and rational cohomology for $G$ identifies with cohomology for the Hopf algebra $kG$.

For $G = \Gar$, one has $k[\Gar] = k[T]/(T^{p^r})$, so
	\begin{equation} \label{eq:kGar}
	k\Gar \cong k[u_0,\ldots,u_{r-1}]/(u_0^p,\ldots,u_{r-1}^p).
	\end{equation}
Here $u_i$ is the functional such that $u_i(f)$ is equal to the coefficient of $T^{p^i}$, for each $f \in k[T]/(T^{p^r})$. Given an integer $0 \leq j < p^r$, let $j = \sum_{\ell=0}^{r-1} j_\ell p^\ell$ be its base-$p$ decomposition (so $0 \leq j_\ell < p$), and set
	\[
	\gamma_j = \frac{u_0^{j_0} \cdot u_1^{j_1} \cdots u_{r-1}^{j_{r-1}}}{j_0! \cdot j_1 ! \cdots j_{r-1}!} \in k\Gar.
	\]
Then the $\gamma_j$ form a $k$-basis for $k\Gar$. The coproduct $\Delta$ and antipode $S$ on $k\Gar$ are given by
	\begin{equation} \label{eq:kGar-coproduct-antipode}
	\Delta(\gamma_\ell) = \sum_{i+j=\ell} \gamma_i \otimes \gamma_j \quad \text{and} \quad S(\gamma_j) = (-1)^j \gamma_j.
	\end{equation}

Given a point $\fs \in V_r(G)$ (i.e., a prime ideal $\fs \in \Spec k[V_r(G)]$), let $k(\fs)$ be the residue field of $V_r(G)$ at $\fs$, and let $\phi_{\fs} \in \Hom_{\CAlg(k)}(k[V_r(G)], k(\fs))$ be the canonical $k$-algebra homomorphism. Then $\phi_{\fs}$ defines a $k(\fs)$-point of $V_r(G)$, and hence determines a homomorphism of $k(\fs)$-group schemes $\nu_{\fs}: \Gar \otimes_k k(\fs) \to G \otimes_k k(\fs)$. Now given a rational $G$-module $M$, let $\nu_{\fs}^*(M \otimes_k k(\fs))$ be the rational $\Gar \otimes_k k(\fs)$-module (equivalently, the $k(\fs)\Gar$-module) obtained via pullback along $\nu_{\fs}$. Then the subset $V_r(G)_M$ of $V_r(G)$ is defined by
	\begin{equation} \label{eq:VrGM} \begin{split}
	V_r(G)_M = \{ &\fs \in V_r(G) : \nu_{\fs}^*(M \otimes_k k(\fs)) \text{ is not projective over the subalgebra} \\
	& k(\fs)[u_{r-1}]/(u_{r-1}^p) \subset k(\fs)[u_0,\ldots,u_{r-1}]/(u_0^p,\ldots,u_{r-1}^p) = k(\fs) \Gar \}.
	\end{split}
	\end{equation}
By \cite[Proposition 6.1]{Suslin:1997a}, if $G$ is an infinitesimal $k$-group scheme and if $M$ is a finite-dimensional rational $G$-module, then $V_r(G)_M$ is a Zariski closed conical subset of $V_r(G)$. If $r = 1$, then the set of $k$-points of $V_r(G)_M$ identifies with the `rank variety' defined by the right-hand side of \eqref{eq:restricted-support-equals-rank}.

The identity map on $k[V_r(G)]$ defines a $k[V_r(G)]$-point of $V_r(G)$, and hence determines a universal homomorphism $u: \Gar \otimes_k k[V_r(G)] \to G \otimes_k k[V_r(G)]$ of group schemes over $k[V_r(G)]$. In \cite{Suslin:1997}, Suslin, Friedlander, and Bendel (SFB) used this map to define, for each affine (algebraic) $k$-group scheme $G$, a natural $k$-algebra homomorphism
	\begin{equation} \label{eq:psir}
	\psi_r: H(G,k) \to k[V_r(G)]
	\end{equation}
Here $H(G,k) = \opH^{\ev}(G,k)$ if $p$ is odd, and $H(G,k) = \Hbul(G,k)$ if $p=2$. Using the algebra relations among the extension classes \eqref{eq:universal-classes}, they also showed that $\phi_{\GLnr}$ factors through a map
	\begin{equation} \label{eq:phi-bar-GLnr}
	\ol{\phi}: k[V_r(\GLnr)] \to H(\GLnr,k).
	\end{equation}

Let $\Psi$ and $\Phi$ denote the morphisms of affine schemes induced by \eqref{eq:psir} and \eqref{eq:phi-bar-GLnr}, respectively. Through a detailed analysis of how the universal extension classes \eqref{eq:universal-classes} restrict to the Frobenius kernels of $\Ga$, SFB determined that the composite morphism
	\[
	\Theta: V_r(\GLnr) \xrightarrow{\Psi} \Spec H(\GLnr,k) \xrightarrow{\Phi} V_r(\GLnr)
	\]
equals the $r$-th Frobenius twist morphism for the scheme $V_r(\GLnr) = V_r(\GLn)$ \cite[Theorem 5.2]{Suslin:1997}. In particular, the homomorphism $\ol{\phi}: k[V_r(\GLnr)] \to H(\GLnr,k)$ has nilpotent kernel, and by naturality it follows that the map $\psi_r: H(G,k) \to k[V_r(G)]$ is surjective onto $p^r$-th powers for any infinitesimal $k$-group scheme $G$ of height $\leq r$.

Next, SFB proved the following detection theorem for infinitesimal group schemes:

\begin{theorem}[{\cite[Theorem 4.3]{Suslin:1997a}}] \label{thm:SFB-thm-4.3}
Let $G$ be an infinitesimal group scheme of height $\leq r$ over $k$ and let $\Lambda$ be an associative unital rational $G$-algebra. Then the following conditions on a cohomology class $z \in \opH^n(G,\Lambda)$ are equivalent:
	\begin{enumerate}
	\item $z$ is nilpotent.
	\item For every field extension $K/k$ and every $K$-group scheme homomorphism $\nu: \Gar \otimes_k K \to G \otimes_k K$, the cohomology class $\nu^*(z_K) \in \opH^n(\Gar \otimes_k K, \Lambda_K)$ is nilpotent.
	\item For every point $\fs \in V_r(G)$, the class $\nu_{\fs}^*(z_{k(\fs)}) \in \opH^n(\Gar \otimes_k k(\fs), \Lambda_{k(\fs)})$ is nilpotent.
	\end{enumerate}
\end{theorem}

The detection theorem was proved first for unipotent infinitesimal group schemes. The general case was then handled with the help of a spectral sequence relating the cohomology of a Frobenius kernel $G_{(r)}$ of a reductive group $G$ to the cohomology of the Frobenius kernel $B_{(r)}$ of a Borel subgroup $B$ of $G$. Similar ideas were also used in the proof of the equality \eqref{eq:restricted-support-equals-rank}.

Using the detection theorem, SFB deduced (for $G$ infinitesimal of height $\leq r$) that the kernel of the $k$-algebra map $\psi_r: H(G,k) \to k[V_r(G)]$ is also nilpotent, and hence $\psi_r$ induces a finite universal homeomorphism of schemes
	\[
	\Psi: V_r(G) \simeq \abs{G} := \Spec H(G,k).
	\]
Furthermore, they deduced for each finite-dimensional rational $G$-module $M$ that $\Psi$ restricts to a homeomorphism
	\[
	\Psi: V_r(G)_M \simeq \abs{G}_M := \Spec\big( H(G,k)/I_G(M) \big).
	\]
The rank variety description for $\abs{G}_M$ afforded by $V_r(G)_M$ implies the Tensor Product Property for finite-dimensional rational $G$-modules $M$ and $N$:
	\[
	\abs{G}_{M \otimes N} = \abs{G}_M \cap \abs{G}_N.
	\]

\subsection{Finite group schemes} \label{subsec:finite-group-schemes}

Let $G$ be a finite $k$-group scheme. A \emph{$\pi$-point} of $G$ is a flat map of $K$-algebras $\alpha: K[t]/(t^p) \to KG_K$, for some field extension $K$ of $k$, which factors through the group algebra $KE$ of an abelian unipotent subgroup scheme $E$ of $G_K$. Two $\pi$-points $\alpha: K[t]/(t^p) \to KG_K$ and $\beta: L[t]/(t^p) \to LG_L$ are equivalent, denoted $\alpha \sim \beta$, if for any finite-dimensional $kG$-module $M$, the pullback $\alpha^*(M_K)$ is projective over $K[t]/(t^p)$ if and only if the pullback $\beta^*(M_L)$ is projective over $L[t]/(t^p)$. Denote the equivalence class of a $\pi$-point $\alpha$ by $[\alpha]$, and denote the set of equivalence classes of $\pi$-points in $G$ by $\Pi(G)$.

Given a finite-dimensional rational $G$-module $M$, set
	\[
	\Pi(G)_M = \set{ [\alpha: K[t]/(t^p) \to KG_K] \in \Pi(G) : \alpha^*(M_K) \text{ is not projective}}.
	\]
Then the class of subsets $\set{ \Pi(G)_M : \text{$M$ finite-dimensional $kG$-module}}$ is the class of closed subsets of a Noetherian topology on $\Pi(G)$.

Given a $\pi$-point $\alpha: K[t]/(t^p) \to KG$, let $\Hbul(\alpha)$ be the composite morphism
	\[
	\Hbul(G,k) = \Ext_G^\bullet(k,k) \xrightarrow{K \otimes_k -} \Ext_{G_K}^\bullet(K,K) \xrightarrow{\alpha^*} \Ext_{K[t]/(t^p)}^\bullet(K,K).
	\]
Then $\sqrt{\ker(\Hbul(\alpha))}$ is a homogeneous prime ideal in $\Hbul(G,k)$, and hence defines a point $\fp(\alpha)$ in the projective spectrum $\Proj \Hbul(G,k)$. Conversely, each point $\fp \in \Proj \Hbul(G,k)$ can be realized as $\fp = \fp(\alpha_{\fp})$ for some $\pi$-point $\alpha_{\fp}: K[t]/(t^p) \to KG_K$ of $G$. In this way, one gets a homeomorphism of topological spaces
	\[
	\Psi_G : \Pi(G) \simeq \Proj \Hbul(G,k), \quad [\alpha] \mapsto \fp(\alpha),
	\]
which satisfies
	\[
	\Psi_G(\Pi(G)_M) = \Proj( \abs{G}_M ) = \Proj\big( \Hbul(G,k) / I_G(M) \big)
	\]
for each finite-dimensional $kG$-module $M$ \cite[Theorem 3.6]{Friedlander:2007}.

For an arbitrary (i.e., not necessarily finite-dimensional) $kG$-module $M$, the $\pi$-support of $M$, denoted $\pisupp_G(M)$, is the subset of $\Proj \Hbul(G,k)$ defined by
	\[
	\pisupp_G(M) = \set{ \fp \in \Proj \Hbul(G,k) : \alpha_{\fp}^*(M_K) \text{ is not projective}}.
	\]
If $M$ is finite-dimensional, then $\pisupp_G(M) = \Psi_G(\Pi(G)_M)$. For arbitrary $kG$-modules, one can show (see \cite{Friedlander:2007,Benson:2022}) that $\pi$-support detects projectivity,
	\[
	\pisupp_G(M) = \emptyset \quad \text{if and only if} \quad \text{$M$ is projective},
	\]
and $\pi$-support satisfies the Tensor Product Property:
	\[
	\pisupp_G(M \otimes N) = \pisupp_G(M) \cap \pisupp_G(N).
	\]

Let $\iota: k[t]/(t^p) \to k\Gar$ be the $k$-algebra map defined by $\iota(t) = u_{r-1}$. If $\nu : \Gar \otimes_k K \to G \otimes_k K$ is a homomorphism defined over a field extension $K/k$, then the induced $K$-algebra map $\nu \circ \iota: K[t]/(t^p) \to K\Gar \to KG_K$ is a $\pi$-point of $G$. In this way, one gets for each infinitesimal $k$-group scheme $G$ of height $\leq r$ a bijection $\Proj V_r(G) \simeq \Pi(G)$, which restricts for each finite-dimensional $kG$-module $M$ to a bijection $\Proj V_r(G)_M \simeq \Pi(G)_M$; see \cite[Proposition 3.8]{Friedlander:2005}.

\section{Infinitesimal super\-group schemes} \label{sec:infinitesimal-supergroups}

In this section, let $r \geq 1$ be a fixed integer.

\subsection{Cohomology of infinitesimal super\-groups}

The general linear Lie super\-algebra $\glmn$ identifies with the space $\Hom_k(k^{m|n},k^{m|n})$ of $k$-linear endomorphisms of the superspace $k^{m|n}$. Then $\glmn_{\zero} = \gl_m \oplus \gl_n$ and $\glmn_{\one} = \g_{+1} \oplus \g_{-1}$, where $\g_{+1} = \Hom_k(k^{0|n},k^{m|0})$ and $\g_{-1} = \Hom_k(k^{m|0},k^{0|n})$. The adjoint action of the general linear super\-group $\GLmn$ restricts to a rational action of $GL_m \times GL_n$ on each of $\gl_m$, $\gl_n$, $\g_{+1}$, and $\g_{-1}$, and the $r$-th Frobenius morphism on $\GLmn$ defines a super\-group scheme homomorphism $F^r : \GLmn \to (GL_m \times GL_n)^{(r)}$. Then via this map, the Frobenius twists $\gl_m^{(r)}$, $\gl_n^{(r)}$, $\g_{+1}^{(r)}$, and $\g_{-1}^{(r)}$ each become rational $\GLmn$-supermodules.

In \cite{Drupieski:2016}, the first author used the theory of strict polynomial superfunctors to establish the existence of certain universal extension classes for the general linear super\-group:
	\begin{equation} \label{eq:super-universal-classes}
	\left.
	\begin{aligned}
	\bse_i^{(r-i)} &\in \opH^{2p^{i-1}}(\GLmn,\gl_m^{(r)}) \\
	(\bse_i^{(r-i)})^\Pi &\in \opH^{2p^{i-1}}(\GLmn,\gl_n^{(r)})
	\end{aligned}
	\right\} \text{ for $1 \leq i \leq r$, and }
	\left\{
	\begin{aligned}
	\bsc_r &\in \opH^{p^r}(\GLmn,\g_{+1}^{(r)}), \\
	\bsc_r^\Pi &\in \opH^{p^r}(\GLmn,\g_{-1}^{(r)}).
	\end{aligned}
	\right.
	\end{equation}
In a manner similar to \eqref{eq:phi-GLnr}, the restriction of these classes to the $r$-th Frobenius kernel of $\GLmn$ can be used to construct a homomorphism of graded super\-algebras
	\begin{equation} \label{eq:phi-GLmnr} \textstyle
	\phi_{\GLmnr}: \left( \bigotimes_{i=1}^r S(\glmn_{\zero}^*[2p^{i-1}])^{(r)} \right) \otimes S( \glmn_{\one}^*[p^r])^{(r)} \to \Hbul(\GLmnr,k).
	\end{equation}
Given a closed sub-super\-group scheme $G$ of $\GLmnr$, one can compose with the restriction map in cohomology $\Hbul(\GLmnr,k) \to \Hbul(G,k)$ to get a homomorphism of graded super\-algebras
	\[ \textstyle
	\phi_G: \left( \bigotimes_{i=1}^r S(\glmn_{\zero}^*[2p^{i-1}])^{(r)} \right) \otimes S( \glmn_{\one}^*[p^r])^{(r)} \to \Hbul(G,k).
	\]
One of the main consequences of the work in \cite{Drupieski:2016} was that, if $G$ is an infinitesimal $k$-super\-group scheme of height $\leq r$, embedded into some $\GLmnr$, then $\phi_G$ is a finite map and $\Hbul(G,M)$ is finite over the image of $\phi_G$ whenever $M$ is a finite-dimensional rational $G$-supermodule. In particular, $\Hbul(G,k)$ is a finitely-generated $k$-super\-algebra, and $\Hbul(G,M)$ is a finite $\Hbul(G,k)$-module.

In the case $r=1$, the map $\phi_G$ can be interpreted via the edge maps of a spectral sequence. Indeed, let $G$ be a height-$1$ infinitesimal $k$-super\-group scheme, and let $\g$ be its restricted Lie super\-algebra. Arguing in exactly the same manner as for \eqref{eq:FPspecseq}, one gets a spectral sequence
	\begin{equation} \label{eq:super-FPspecseq}
	E_2^{i,j} = S^{i/2}(\gzero^*)^{(1)} \otimes \opH^j(\g,k) \Rightarrow \opH^{i+j}(G,k).
	\end{equation}
Fixing a choice of embedding $G \hookrightarrow GL_{m|n(1)}$, one gets a map of restricted Lie super\-algebras $\g \hookrightarrow \glmn$. Then up to a nonzero scalar, the restricted map $\phi_G: S(\glmn_{\zero}^*)^{(1)} \to \Hbul(G,k)$ is the composition of the quotient map $S(\glmn_{\zero}^*)^{(1)} \twoheadrightarrow S(\gzero^*)^{(1)}$ and the horizontal edge map of \eqref{eq:super-FPspecseq}. Similarly, the composition of $\phi_G: S(\glmn_{\one}^*)^{(1)} \to \Hbul(G,k)$ with the horizontal edge map $\Hbul(G,k) \to \Hbul(\g,k)$ of \eqref{eq:super-FPspecseq} is equal, up to a nonzero scalar, to the composition of the quotient map $S(\glmn_{\one}^*)^{(1)} \twoheadrightarrow S(\gone^*)^{(1)}$ and the homomorphism $\varphi: S(\gone^*)^{(1)} \to \Hbul(\g,k)$ of \eqref{eq:varphi}.

Unlike the classical situation of \eqref{eq:FPspecseq}, in general the $E_2$-page of \eqref{eq:super-FPspecseq} is not finite over the evident subalgebra of permanent cycles generated by the terms on the horizontal axis. This shows how, starting immediately with the base case of $r=1$, the problem of establishing the basic finite-generation property \ref{FG-ring} is more subtle for infinitesimal supergroups than for ordinary infinitesimal group schemes.

\subsection{Multiparameter super\-groups} \label{subsec:multiparameter}

Define $\Mr$ to be the affine $k$-super\-group scheme whose coordinate algebra $k[\Mr]$ is generated by the odd element $\tau$ and the even elements $\theta$ and $\sigma_i$ for $i \in \N$, such that $\tau^2 = 0$, $\sigma_0 = 1$, $\theta^{p^{r-1}} = \sigma_1$, and $\sigma_i \sigma_j = \binom{i+j}{i} \sigma_{i+j}$, i.e.,
	\[
	k[\Mr] = k[\tau,\theta,\sigma_1,\sigma_2,\sigma_3,\ldots]/\big( \tau^2, \theta^{p^{r-1}}-\sigma_1,\sigma_i\sigma_j - \tbinom{i+j}{i}\sigma_{i+j}: i,j \in \N \big).
	\]
Then the set of monomials $\set{ \theta^i \sigma_j, \tau \theta^i \sigma_j : 0 \leq i < p^{r-1}, j \in \N}$ is a homogeneous basis for $k[\Mr]$. The coproduct $\Delta$ and the antipode $S$ on $k[\Mr]$ are defined on generators by the formulas
	\begin{align*}
	\Delta(\tau) &= \tau \otimes 1 + 1 \otimes \tau, & S(\tau) &= -\tau, \\
	\Delta(\theta) &= \theta \otimes 1 + 1 \otimes \theta, & S(\theta) &= -\theta, \\
	\Delta(\sigma_i) &= \textstyle \sum_{u+v=i} \sigma_u \otimes \sigma_v + \sum_{u+v+p=i} \sigma_u \tau \otimes \sigma_v \tau, & S(\sigma_i) &= (-1)^i \sigma_i.
	\end{align*}
The $\Z_2$-grading on $k[\Mr]$ lifts to a $\Z$-grading such that $\deg(\tau) = p^r$, $\deg(\theta) = 2$, and $\deg(\sigma_i) = 2ip^{r-1}$. This makes $k[\Mr]$ into a $\Z$-graded Hopf algebra of finite type.

Let $k\Mr = k[\Mr]^*$ be the `group algebra' of $\Mr$, i.e., the $k$-algebra that is dual to the coalgebra structure on $k[\Mr]$. For $0 \leq i \leq r-1$, let $u_i \in k\Mr$ be the (even) linear functional that is dual to the basis vector $\theta^{p^i}$ (so in particular, $u_{r-1}$ is dual to $\theta^{p^{r-1}} = \sigma_1$), and let $v \in k\Mr$ be the (odd) linear functional that is dual to the basis vector $\tau$. Then
	\[
	k\Mr = \frac{k[u_0,\ldots,u_{r-2},v][[u_{r-1}]]}{(u_0^p,\ldots,u_{r-2}^p,u_{r-1}^p+v^2)}.
	\]
The graded dual of the $\Z$-graded vector space $k[\Mr]$ is the subalgebra
	\[
	\Pr = \frac{k[u_0,\ldots,u_{r-1},v]}{(u_0^p,\ldots,u_{r-2}^p,u_{r-1}^p+v^2)}.
	\]
The algebra $\Pr$ can be thought of as the `polynomial part' of $k\Mr$. Since $k[\Mr]$ is a $\Z$-graded Hopf algebra of finite type, its graded dual is also a $\Z$-graded Hopf algebra (and hence a Hopf super\-algebra, by reducing the $\Z$-grading modulo $2$). The coproduct $\Delta$ and antipode $S$ for $\Pr$ are defined on the generator $v$ by $\Delta(v) = v \otimes 1 + 1 \otimes v$ and $S(v) = -v$, and are defined on $u_0,\ldots,u_{r-1}$ by the same formulas as in \eqref{eq:kGar-coproduct-antipode}. In particular, $\Delta(u_{r-1}^p) = u_{r-1}^p \otimes 1 + 1 \otimes u_{r-1}^p$.

Let $f = \sum_{i=1}^t c_i T^{p^i}$ be an inseparable $p$-polynomial (i.e., a $p$-polynomial without a linear term), and let $\eta \in k$. Since $u_0$ and $u_{r-1}$ are both primitive in $\Pr$, the sum $f(u_{r-1})+ \eta \cdot u_0$ is also primitive in $\Pr$. Assuming that $f \neq 0$ if $r \geq 2$, and that either $f \neq 0$ or $\eta \neq 0$ if $r = 1$, the quotient
	\[
	k\Mrfeta := \Pr/( f(u_{r-1})+\eta \cdot u_0)
	\]
is then a finite-dimensional Hopf super\-algebra. Let $\Mrfeta$ be the finite $k$-super\-group scheme whose group algebra is $k\Mrfeta$. For $s \geq 1$, set $\Mrseta = \M_{r,T^{p^s},\eta}$, and set $\Mrs = \M_{r,s,0}$. Then
	\begin{align*}
	k\Mrseta &= \frac{k[u_0,\ldots,u_{r-1},v]}{(u_0^p,\ldots,u_{r-2}^p,u_{r-1}^p+v^2,u_{r-1}^{p^s}+\eta \cdot u_0 )}, \quad \text{and} \\
	k\Mrs &= \frac{k[u_0,\ldots,u_{r-1},v]}{(u_0^p,\ldots,u_{r-2}^p,u_{r-1}^p+v^2,u_{r-1}^{p^s} )}.
	\end{align*}
Observe also that $\P_{r-1} \cong \Pr/(u_0)$, $k\Gar \cong \Pr/(v)$, and
	\[
	k\Gam = k[v]/(v^2) \cong \Pr/(u_0,\ldots,u_{r-1}).
	\]
Here $\Gam$ is the odd additive super\-group scheme. Its coordinate super\-algebra is $k[\Gam] = k[t]/(t^2)$ with $\ol{t} = \one$. As a functor $\Gam: \CSAlg(k) \to \Grp$, it is defined by $\Gam(A) = (\Aone,+)$. Here $\CSAlg(k)$ denotes the category of (super)commutative $k$-super\-algebras.

We say that a finite $k$-super\-group scheme $G$ is a \emph{multiparameter super\-group of height $\leq r$} if its group algebra $kG = k[G]^*$ is isomorphic to a (finite-dimensional) Hopf super\-algebra quotient of $\Pr$. By \cite[Proposition 2.2.1]{Drupieski:2021a}, each multiparameter super\-group of height $\leq r$ is isomorphic to (at least) one of the following:
	\begin{enumerate}[wide=\parindent,label={(MP\arabic*)}]
	\item $\Gas$ for some $0 \leq s \leq r$,
	\item $\Gas \times \Gam$ for some $0 \leq s \leq r$ (note that $\Gas \times \Gam = \M_{s,1}$), or
	\item $\Msfeta$ for some $1 \leq s \leq r$, some inseparable $p$-polynomial $0 \neq f \in k[T]$, and some $\eta \in k$.
	\end{enumerate}
Among these supergroups, the \emph{unipotent} multiparameter supergroups have the forms:
	\begin{enumerate}[wide=\parindent,label={(UMP\arabic*)}]
	\item \label{ump-Gar} $\Gar$ for some $r \geq 0$,
	\item \label{ump-Gar-Gam} $\Gar \times \Gam$ for some $r \geq 0$,
	\item \label{ump-Mrs} $\Mrs$ for some $r,s \geq 1$, or
	\item \label{ump-Mrseta} $\Mrseta$ for some integers $r \geq 2$, $s \geq 1$, and some scalar $0 \neq \eta \in k$.
	\end{enumerate}
By definition, $\G_{a(0)}$ is the trivial (constant) group scheme. 
	
\subsection{The functor of multiparameter super\-groups}

Given affine $k$-super\-group schemes $G$ and $H$, let $\Hom(G,H): \CAlg(k) \to \Set$ be the functor defined by
	\[
	\Hom(G,H)(A) = \Hom_{\sGrp/A}(G \otimes_k A, H \otimes_k A),
	\]
the set of $A$-super\-group scheme homomorphisms $\rho: G \otimes_k A \to H \otimes_k A$. Similarly, given Hopf $k$-super\-algebras $R$ and $S$, let $\Hom(R,S) : \CAlg(k) \to \Set$ be the functor defined by
	\[
	\Hom(R,S)(A) = \Hom_{\sHopf/A}(R \otimes_k A, S \otimes_k A),
	\]
the set of Hopf $A$-super\-algebra homomorphisms $\rho: R \otimes_k A \to S \otimes_k A$. If $G$ and $H$ are finite $k$-super\-group schemes, then $\Hom(G,H) \cong \Hom(kG,kH)$. 

\begin{definition} \label{def:VrG-super}
Given a finite $k$-super\-group scheme $G$, define $V_r(G) : \CAlg(k) \to \Set$ by
	\[
	V_r(G) = \Hom(\Pr,kG).
	\]
We call $V_r(G)$ the \emph{functor of multiparameter super\-groups of height $\leq r$ in $G$}.
\end{definition}

We consider the objects in $\CAlg(k)$ as superalgebras concentrated in even superdegree. If $G$ is an ordinary (purely even) finite $k$-group scheme, then any $\rho \in V_r(G)(A)$ will automatically factor through the canonical quotient map $\Pr \otimes_k A \twoheadrightarrow \Pr/(v) \otimes_k A \cong k\Gar \otimes_k A$, and hence will define a map of $A$-group schemes $\rho: \Gar \otimes_k A \to G \otimes_k A$ (and conversely). Thus when $G$ is purely even, our new use of the notation $V_r(G)$ agrees with that already established in Section \ref{sec:SFB-support-varieties} (although our new usage requires $G$ to be finite, whereas before $G$ was merely assumed to be affine algebraic, to get a scheme structure on $V_r(G)$).

By \cite[Lemma 4.1.1]{Drupieski:2021a}, the functor $V_r(G)$ admits the structure of an affine $k$-scheme of finite type, and the assignment $G \mapsto V_r(G)$ is then a covariant functor from finite $k$-super\-group schemes to affine $k$-schemes of finite type, which takes closed embeddings (of super\-group schemes) to closed embeddings.\footnote{In \cite{Drupieski:2021a}, we defined a $k$-\emph{super}functor, denoted $\bsvr(G)$, which admitted the structure of an affine $k$-\emph{super}scheme of finite type, and whose purely underlying even subscheme is $V_r(G)$. The superfunctor and superscheme structures can be ignored as far as the applications to cohomological support varieties are concerned, so we dispense with them in the current article.} If $E$ is a multiparameter supergroup of height $\leq r$ (i.e., if $kE$ is a Hopf super\-algebra quotient of $\Pr$), then the functor $\Hom(kE,kG) \cong \Hom(E,G)$ also admits an affine $k$-scheme structure, and identifies via the quotient map $\Pr \twoheadrightarrow kE$ with a closed subscheme of $V_r(G)$.\footnote{More generally, if $G$ is an affine algebraic $k$-group scheme, i.e., if $k[G]$ is finitely-generated as a $k$-algebra but not necessarily finite-dimensional, and if $E$ is a multiparameter supergroup, then $\Hom(E,G)$ admits an affine $k$-scheme structure by \cite[Lemma 3.3.6]{Drupieski:2019b}.} Identifying $\Hom(E,G)$ with its image in $V_r(G)$ via this embedding, one then has
	\begin{equation} \label{eq:VrG-union-homs}
	V_r(G)(k) = \bigcup_{(E \leq G)/\cong} \Hom(E,G)(k),
	\end{equation}
where the union is over all isomorphism classes of multiparameter $k$-supergroups of height $\leq r$ that occur as closed subsupergroups of $G$. Similarly one has
	\begin{equation} \label{eq:VrG-union-VrE}
	V_r(G)(k) = \bigcup_{E \leq G} V_r(E)(k),
	\end{equation}
where the union is over all multiparameter subsupergroups of height $\leq r$ in $G$  \cite[Lemma 4.1.4]{Drupieski:2021a}.

\begin{remark} \label{rem:Hom(Mr,G)}
For all $1 \leq s \leq s'$, there are canonical quotient maps $\Mr \twoheadrightarrow \M_{r,s'} \twoheadrightarrow \M_{r,s}$. If $G$ is an algebraic $k$-super\-group scheme and if $\rho \in \Hom(\Mr,G)(A)$, then $\rho$ must factor through the quotient $\Mr \otimes_k A \twoheadrightarrow \Mrs \otimes_k A$ for some $s \geq 1$ \cite[Remark 3.1.3(4)]{Drupieski:2019a}. In this way one gets
	\begin{equation} \label{eq:Hom(Mr,G)}
	\Hom(\Mr,G) = \bigcup_{s \geq 1} \Hom(\Mrs,G).
	\end{equation}
If $G$ is a finite unipotent $k$-supergroup scheme, this union is finite; see Section \ref{subsec:spectrum-unipotent}. In the papers \cite{Drupieski:2019a,Drupieski:2019c} we initially considered the functor $\Hom(\Mr,G)$ as a possible superization of the classical functor of one-parameter subgroups \eqref{eq:Vrg-non-super}, before deciding while writing \cite{Drupieski:2021a} that $\Hom(\Pr,kG)$ seems more likely to provide a better generalization for arbitrary finite supergroups, at least as far as applications to support varieties are concerned.
\end{remark}

\begin{question}
Is there an affine (though perhaps not algebraic) $k$-supergroup scheme $\G$ such that for all finite $k$-supergroup schemes $G$, $\Hom(\G,G) \cong \Hom(\Pr,kG)$?
\end{question}

Given a multiparameter $k$-supergroup $E$ and a finite $k$-supergroup scheme $G$, set
	\[
	V_E(G) = \Hom(E,G).
	\]
Then $V_E(G)$ is a closed subscheme of $V_r(G)$. In particular, set
	\[
	\Vrfeta(G) = \Hom(\Mrfeta,G) \quad \text{and} \quad \Vrs(G) = \Hom(\Mrs,G).
	\]

\subsection{Commuting varieties}

If $G$ is an affine algebraic $k$-super\-group scheme with restricted Lie super\-algebra $\g$, the Lie super\-algebra of $G_A = G \otimes_k A$ is then given by $\g_A = \g \otimes_k A$. Under this identification, the $p$-operation on $(\g_A)_{\zero} = \gzero \otimes_k A$ is given by $(x \otimes a)^{[p]} = x^{[p]} \otimes a^p$.

\begin{definition}
Given an affine $k$-super\-group scheme $G$ with restricted Lie super\-algebra $\g$, let $C_r(G): \CAlg(k) \to \Set$ be the functor defined by
	\begin{multline*}
	C_r(G)(A) = \Big\{ (\alpha_0,\ldots,\alpha_{r-1},\beta) \in [(\g_A)_{\zero}]^{\times r} \times (\g_A)_{\one}: [\alpha_i,\alpha_j] = [\alpha_i,\beta] = 0 \text{ for all $i$ and $j$,} \\
	\alpha_i^{[p]} = 0 \text{ for all $0 \leq i \leq r-2$, and } \alpha_{r-1}^{[p]} + \tfrac{1}{2}[\beta,\beta] = 0 \Big\}.
	\end{multline*}
\end{definition}

\begin{remark} \label{rem:notation-caution}
In the papers \cite{Drupieski:2019b,Drupieski:2019c}, we considered a $k$-superfunctor $\bsvr(\GLmn)$ whose restriction to $\CAlg(k)$ we denoted there by $V_r(\GLmn)$, but which is denoted in this paper by $C_r(\GLmn)$. We advise the reader to exercise caution when comparing the notation in \cite{Drupieski:2019b,Drupieski:2019c} with the sometimes inconsistent notation in \cite{Drupieski:2021a} and in the present paper.
\end{remark}

Since the $p$-operation on $(\glmn_A)_{\zero} \cong \gl_m(A) \oplus \gl_n(A)$ is given simply by raising matrices to the $p$-th power, it readily follows that $C_r(\GLmn)$ is represented by an affine scheme of finite type over $k$, which we also denote $C_r(\GLmn)$. More generally, if $G$ is an affine algebraic $k$-super\-group scheme, one can use a choice of closed embedding $G \hookrightarrow \GLmn$ to show that $C_r(G)$ admits a scheme structure, and that $G \mapsto C_r(G)$ is then a covariant functor from affine algebraic $k$-super\-group schemes to affine $k$-schemes of finite type, which takes closed embeddings (of super\-group schemes) to closed embeddings.

Given an inseparable $p$-polynomial $0 \neq f \in k[T]$ and a scalar $\eta \in k$, one has
	\[
	\Vrfeta(\GLmnr)(A) \cong \big\{ (\alpha_0,\ldots,\alpha_{r-1},\beta) \in C_r(\GLmnr)(A) : f(\alpha_{r-1}) + \eta \cdot \alpha_0 = 0 \big\}
	\]
by \cite[Proposition 3.3.5]{Drupieski:2019a}. In particular, given an integer $s \geq 1$, one has
	\[
	\Vrs(\GLmnr)(A) \cong \big\{ (\alpha_0,\ldots,\alpha_{r-1},\beta) \in C_r(\GLmnr)(A) : \alpha_{r-1}^{p^s} = 0 \big\}.
	\]
For the ordinary (purely even) $k$-group scheme $\GLnr$, one has
	\[
	V_r(\GLnr) = \Hom(\Pr,k\GLnr) \cong \Hom(\Gar,\GLnr) \cong C_r(\GLnr)
	\]
by \eqref{eq:Vr(GLn)}. More generally, if $G$ is a finite purely even $k$-group scheme that admits an embedding $G \hookrightarrow GL_n$ of \emph{exponential type} \cite[p.\ 697]{Suslin:1997}, then $V_r(G) \cong C_r(G)$ by \cite[Lemma 1.7]{Suslin:1997}.

\begin{remark} \label{rem:V1G-height-1}
Let $G$ be a height-$1$ infinitesimal $k$-supergroup scheme, and let $\g$ be its (restricted) Lie superalgebra. Then $kG = V(\g)$. In this situation, a Hopf super\-algebra homomorphism $\sigma: \Pone = k[u,v]/(u^p+v^2) \to kG = V(\g)$ is determined by the data of the even primitive element $\sigma(u) \in \gzero$ and the odd primitive element $\sigma(v) \in \gone$, which must satisfy $\sigma(u)^p + \sigma(v)^2 = 0$. Then
	\[
	V_1(G)(k) = \Hom_{\sHopf/k}(\Pone,V(\g)) \cong \big\{ (\alpha,\beta) \in \gzero \times \gone : \alpha^{[p]} + \tfrac{1}{2}[\beta,\beta] = 0 \big\}.
	\]
In particular, $V_1(G)(k) \cong C_1(G)(k)$ for $G$ of height $1$.
\end{remark}

\subsection{The cohomological spectrum}

\subsubsection{The general case}

Let $G$ be a finite $k$-supergroup scheme. The identity map on $k[V_r(G)]$ defines a $k[V_r(G)]$-point of $V_r(G)$, and hence determines a universal homomorphism $u_G: \Pr \otimes_k k[V_r(G)] \to kG \otimes_k k[V_r(G)]$ of Hopf superalgebras over $k[V_r(G)]$. In \cite[\S4.2]{Drupieski:2021a}, we used this map to define a natural (with respect to $G$) $k$-algebra homomorphism
	\begin{equation} \label{eq:psir-super}
	\psi_r: H(G,k) \to k[V_r(G)].
	\end{equation}
The ring $k[V_r(G)]$ is naturally $\Z[\frac{p^r}{2}]$-graded, and the map $\psi_r$ then multiplies $\Z$-degrees by $\frac{p^r}{2}$. Composing \eqref{eq:psir-super} with the quotient map $k[V_r(G)] \twoheadrightarrow k[V_E(G)]$ for a multiparameter $k$-supergroup $E$, one gets a $k$-algebra homomorphism
	\begin{equation} \label{eq:psirfeta}
	\psi_E : H(G,k) \to k[V_E(G)].
	\end{equation}
Set $\psi_{r,f,\eta} = \psi_{\Mrfeta}$, and set $\psi_{r,s} = \psi_{\Mrs} = \psi_{r,T^{p^s},0}$.

The domain of the super\-algebra homomorphism \eqref{eq:phi-GLmnr} identifies with the coordinate algebra of the affine space $(\bigoplus_{i=1}^r \glmn_{\zero}) \oplus \glmn_{\one}$. In \cite[Proposition 6.1.1]{Drupieski:2019b}, we used the algebra relations among the extension classes \eqref{eq:super-universal-classes} to show that, analogously to \eqref{eq:phi-bar-GLnr}, the map $\phi_{\GLmnr}$ factors through a super\-algebra homomorphism
	\begin{equation} \label{eq:phi-bar-GLmnr}
	\ol{\phi}_{\GLmnr}: k[C_r(\GLmnr)] \to H(\GLmnr,k).
	\end{equation}

Let $\Psi_E$ and $\Phi$ denote the morphisms of affine schemes induced by \eqref{eq:psirfeta} and \eqref{eq:phi-bar-GLmnr}, respectively. Through an analysis of how the universal extension classes \eqref{eq:super-universal-classes} restrict to multi\-parameter super\-groups, we showed in \cite[Theorem 6.2.3]{Drupieski:2019b} that the composite morphism
	\begin{equation} \label{eq:theta-morphism}
	\Theta_{r,f,\eta}: \Vrfeta(\GLmnr) \xrightarrow{\Psi_{r,f,\eta}} \Spec H(\GLmnr,k) \xrightarrow{\Phi} C_r(\GLmnr)
	\end{equation}
is equal to the composite of the closed embedding $\Vrfeta(\GLmnr) \hookrightarrow V_r(\GLmn)$ and the $r$-th Frobenius twist morphism for the scheme $V_r(\GLmnr) = V_r(\GLmn)$. In particular, the map $\psi_{r,f,\eta}$ is surjective onto $p^r$-th powers, and the kernel of $\ol{\phi}$ is nilpotent. Using naturality of $\psi_{r,f,\eta}$ with respect to $G$, we deduce for an arbitrary infinitesimal $k$-supergroup scheme $G$ of height $\leq r$ that the map $\psi_{r,f,\eta}: H(G,k) \to k[V_{r,f,\eta}(G)]$ is surjective onto $p^r$-th powers.

\begin{remark}
Allowing $f$ to vary, we used the morphisms $\Theta_{r,f,0}$ to deduce in \cite[Corollary 6.2.4]{Drupieski:2019b} that the finite morphism of affine algebraic varieties
	\[
	\Phi: \MaxSpec H(\GLmnr,k) \to C_r(\GLmnr)(k) = C_r(\GLmn)(k)
	\]
induced by $\Phi$ is a surjection. Now let $G$ be a closed subsupergroup of $\GLmnr$. Composing \eqref{eq:phi-bar-GLmnr} with the restriction map in cohomology, one gets a $k$-algebra map $\ol{\phi}_G: k[C_r(\GLmnr)] \to H(G,k)$. In the special case $r=1$, we showed in \cite[\S 5.4]{Drupieski:2019a} that $\ol{\phi}_G$ factors through $k[C_1(G)]$, and that the induced finite morphism of varieties
	\[
	\Phi_G: \MaxSpec H(G,k) \to C_1(G)(k) = \big\{ (\alpha,\beta) \in \gzero \times \gone : [\alpha,\beta] = 0 \text{ and } \alpha^{[p]} + \tfrac{1}{2}[\beta,\beta] = 0 \big\}
	\]
is a surjection. This provides an analogue for restricted Lie superalgebras of \eqref{eq:image-spectrum-nullcone}.
\end{remark}

\subsubsection{The unipotent case} \label{subsec:spectrum-unipotent}

Now suppose $G$ is an infinitesimal \emph{unipotent} $k$-supergroup scheme of height $\leq r$. Then by \cite[Lemma 4.1.7]{Drupieski:2021a}, there exists an integer $s = s(G) \geq 1$ such that the canonical quotient maps $\Pr \twoheadrightarrow k\M_{r,s'} \twoheadrightarrow k\Mrs$ induce identifications
	\begin{equation} \label{eq:Vrs(G)-stabilize}
	\Hom(\Mrs,G) = \Hom(\M_{r,s'},G) = \Hom(\Pr,kG) = V_r(G) \quad \text{for all $s' \geq s = s(G)$.}
	\end{equation}
So for $s = s(G)$ one has $V_r(G) = \Vrs(G)$, and the homomorphism $\psi_r$ of \eqref{eq:psir-super} is equal to $\psi_{r,s}$, which as we stated above is surjective onto $p^r$-th powers.

Benson, Iyengar, Krause, and Pevtsova (BIKP) \cite{Benson:2021} define a finite supergroup scheme to be \emph{elementary} if it is isomorphic, for some positive integers $r,s,t$, to a quotient of $\Mrs \times (\Z/p)^t$. Here $\Z/p$ denotes the (purely even) constant group scheme corresponding to the cyclic group of order $p$. It turns out that a supergroup scheme is elementary if and only if it is a product of one of the unipotent multiparameter super\-groups listed in \ref{ump-Gar}--\ref{ump-Mrseta} and the constant group scheme $(\Z/p)^t$ for some $t \geq 0$. BIKP then prove the following detection theorem:

\begin{theorem}[{\cite[Theorem 1.2, Theorem 11.2]{Benson:2021}}] \label{thm:BIKP-detection}
Let $G$ be a finite unipotent supergroup scheme over a field $k$ of positive characteristic $p \geq 3$. Then the following hold:
	\begin{enumerate}
	\item \label{item:detect-cohomology-triv-coeff} An element $x \in \Hbul(G,k)$ is nilpotent if and only if, for every extension field $K/k$ and every elementary sub-supergroup scheme $E$ of $G_K$, the restriction of $x_K \in \Hbul(G_K,K)$ to $\Hbul(E,K)$ is nilpotent.
	\item \label{item:detect-projectivity} A $kG$-supermodule $M$ is projective if and only if, for every extension field $K/k$ and every elementary sub-supergroup scheme $E$ of $G_K$, the restriction of $M_K$ to $E$ is projective.
	\item \label{item:detect-cohomology-lambda} Let $\Lambda$ be a unital $G$-algebra. Then an element $x \in \Hbul(G,\Lambda)$ is nilpotent if and only if, for every extension field $K/k$ and every elementary sub-supergroup scheme $E$ of $G_K$, the restriction of $x_K \in \Hbul(G_K, \Lambda_K)$ to $\Hbul(E,\Lambda_K)$ is nilpotent.
	\end{enumerate}
\end{theorem}

Using part \eqref{item:detect-cohomology-triv-coeff} of the theorem, one can then prove:

\begin{theorem}[{\cite[Theorem 5.1.3]{Drupieski:2021a}}] \label{theorem:spectrum-Psir-homeo}
Let $G$ be an infinitesimal unipotent $k$-supergroup scheme of height $\leq r$. Then the kernel of the homomorphism $\psi_r: H(G,k) \to k[V_r(G)]$ is a locally nilpotent ideal, and its image contains all $p^r$-th powers. Consequently, the associated morphism of schemes $\Psi_r$ defines a universal homeomorphism
	\[
	\Psi_r : V_r(G) \simeq \abs{G} = \Spec H(G,k).
	\]
\end{theorem}

\begin{question} \label{question:psir-homeomorphism}
Let $G$ be an arbitrary infinitesimal $k$-supergroup scheme of height $\leq r$. Does the map $\psi_r: H(G,k) \to k[V_r(G)]$ induce a homeomorphism $V_r(G) \simeq \Spec H(G,k)$?
\end{question}

A (or more likely, \emph{the}) key step toward answering Question \ref{question:psir-homeomorphism} would be to establish an analogue of Theorem \ref{thm:BIKP-detection} for non-unipotent $G$. As observed in \cite[Example 1.4.2]{Drupieski:2019b}, the family of elementary supergroups is inadequate for detecting nilpotence of cohomology classes for the non-unipotent supergroup $\M_{1,T^p,-1}$, whose group algebra is $k\M_{1,T^p,-1} = k[u,v]/(u^p+v^2, u^p-u)$.

\begin{question} \label{question:detecting-family}
Does an analogue of Theorem \ref{thm:BIKP-detection} hold for non-unipotent $G$ if the family of elementary supergroups is replaced by the family of supergroups of the form $\M \times (\Z/p)^t$, where $\M$ is one of the multiparameter supergroups defined in Section \ref{subsec:multiparameter}?
\end{question}

As a starting point for trying to answer Question \ref{question:detecting-family}, one might try to imitate in some way the proof of either Theorem \ref{thm:SFB-thm-4.3} or its generalization to arbitrary finite group schemes given in \cite{Suslin:2006}. Roughly, these arguments first establish a detection theorem for unipotent groups, and hence for Frobenius kernels $B_{(r)}$ of the Borel subgroups of $GL_n$. The arguments then exploit algebro-geometric relationships between $GL_n$ and its Borel subgroups to construct a spectral sequence that allows one to bootstrap information from $B_{(r)}$ to all of $GL_{n(r)}$. There are non-trivial difficulties in trying to directly translate this strategy over to $GL_{m|n}$, however, because of such issues as the existence of non-conjugate Borel subgroups, and the lack in general of a super-analogue of the Kempf vanishing theorem. It could be necessary to replace the Borel subgroups in this picture with some other class of subgroups of $GL_{m|n}$ for whom the nilpotence of cohomology classes is detected by restriction to subgroups of the form $\M \times (\Z/p)^t$.

\subsection{Support schemes for supermodules} \label{subsec:support-varieties-super}

Write $\Pone = k[u,v]/(u^p+v^2)$, and let $\iota: \Pone \rightarrow \Pr$ be the superalgebra map defined by $\iota(u) = u_{r-1}$ and $\iota(v) = v$. Let $G$ be a finite $k$-supergroup scheme.  Now given a $kG$-super\-module $M$ and a point $\fs \in V_r(G)$ (i.e., a prime ideal in $\Spec k[V_r(G)]$), we consider $M \otimes_k k(\fs)$ as a $\Pone \otimes_k k(\fs)$-super\-module by pulling back along the composite $k(\fs)$-super\-algebra homomorphism
	\[
	\nu_{\fs} \circ (\iota \otimes 1): \Pone \otimes_k k(\fs) \hookrightarrow \Pr \otimes_k k(\fs) \rightarrow kG \otimes_k k(\fs).
	\]

Recall that the ring $k[V_r(G)]$ is naturally $\Z[\frac{p^r}{2}]$-graded. We say that a Zariski closed subset $X \subseteq V_r(G)$ is \emph{conical} if it is defined by a $\Z[\frac{p^r}{2}]$-homogeneous ideal.

\begin{proposition}[{\cite[Proposition 4.3.1]{Drupieski:2021a}}] \label{prop:VrgMclosed}
Let $G$ be a finite $k$-supergroup scheme, and let $M$ be a $kG$-supermodule. Then
	\[
	\Vrg_M := \big\{ \fs \in \Vrg : \projdim_{\Pone \otimes_k k(\fs)}(M \otimes_k k(\fs)) = \infty \big\}
	\]
is a Zariski closed conical subset of $\Vrg$.
\end{proposition}

Given a $\Pone$-supermodule $M$, the matrices
	\[
	\varphi = \begin{pmatrix} u^{p-1} & v \\ v & -u \end{pmatrix} \quad \text{and} \quad \psi = \begin{pmatrix} u & v \\ v & -u^{p-1} \end{pmatrix}
	\]
naturally define operators $\varphi_M: M \oplus \Pi(M) \to M \oplus \Pi(M)$ and $\psi_M: M \oplus \Pi(M) \to M \oplus \Pi(M)$. Here $\Pi(M)$ denotes the parity shift of $M$. Then for $M$ finite-dimensional, one gets
	\begin{equation} \label{eq:projdim-rank}
	\projdim_{\Pone}(M) = \infty \quad \text{if and only if} \quad \rank(\varphi_M) = \rank(\psi_M) < \dim_k(M);
	\end{equation}
see \cite[\S8]{Benson:2022} or \cite[Remark 4.3.4]{Drupieski:2021a}.

Again making crucial use of the BIKP detection theorem, we proved:

\begin{theorem}[{\cite[Theorem 5.4.1]{Drupieski:2021a}}] \label{thm:Psirinv}
Let $G$ be an infinitesimal unipotent $k$-supergroup scheme of height $\leq r$, and let $M$ be a finite-dimensional rational $G$-supermodule. Then the morphism of schemes $\Psi_r: \Vrg \rightarrow \abs{G}$ satisfies $\Psi_r^{-1}(\abs{G}_M) = \Vrg_M$. Thus, $\Psi_r$ restricts to a homeomorphism
	\begin{equation} \label{eq:VrGM-homeo}
	\Psi_r: \Vrg_M \stackrel{\sim}{\rightarrow} \abs{G}_M.
	\end{equation}
\end{theorem}

\begin{remark} \label{rem:p-nilpotent-explicit}
Combining Remark \ref{rem:V1G-height-1}, the rank criterion \eqref{eq:projdim-rank}, and the theorem, one gets a fairly explicit description of support varieties for height-$1$ infinitesimal unipotent supergroup schemes, or equivalently, for finite-dimensional $p$-nilpotent restricted Lie superalgebras.
\end{remark}

The equalities \eqref{eq:VrG-union-homs} and \eqref{eq:VrG-union-VrE} describe how $\Vrg(k)$ is stratified by pieces coming from the multi\-parameter sub\-super\-group schemes of $G$. The next result, which is reminiscent of Quillen's stratification theorem for finite groups \cite{Quillen:1971}, translates this to $\abs{G}_M$.

\begin{theorem}[{\cite[Theorem 6.1.5]{Drupieski:2021a}}] \label{theorem:supportvarietyunion}
Let $G$ be an infinitesimal unipotent $k$-supergroup scheme of height $\leq r$, and let $M$ be a finite-dimensional rational $G$-supermodule. Then the support variety $\abs{G}_M$ (i.e., the set of $k$-points of the scheme of the same name) can be written as
	\begin{equation} \label{eq:support-stratified}
	\abs{G}_M = \bigcup_{E \leq G} \res_{G,E}^*(\abs{E}_M),
	\end{equation}
where the union is taken over all (closed) elementary $k$-subsupergroup schemes $E$ of $G$, and $\res_{G,E}: H(G,k) \rightarrow H(E,k)$ is the restriction map induced by the embedding $E \hookrightarrow G$.
\end{theorem}

\begin{question}
Do \eqref{eq:VrGM-homeo} and \eqref{eq:support-stratified} remain true as written for non-unipotent $G$?
\end{question}

\subsection{\texorpdfstring{$\pi$-points}{pi-points} for finite supergroup schemes}

In \cite{Benson:2022}, Benson, Iyengar, Krause, and Pevtsova (BIKP) recently introduced a notion of $\pi$-points for arbitrary finite supergroup schemes, generalizing the notion of $\pi$-points for finite group schemes recalled in Section \ref{subsec:finite-group-schemes}.

Set $A_k = k[u,v]/(u^p-v^2)$. Modulo the rescaling $u \mapsto -u$, this is the superalgebra $\Pone$ defined in Section \ref{subsec:multiparameter}. Given a field extension $K/k$, set $A_K = A \otimes_k K$. Given a finite $k$-supergroup scheme $G$, BIKP define a $\pi$-point of $G$ to be a $K$-algebra map of finite flat (equivalently, projective) dimension $\alpha: A_K \to KG_K$, for some field extension $K/k$, such that $\alpha$ factors through the group algebra $KE$ of an elementary sub-supergroup scheme $E$ of $G_K$. Two $\pi$-points $\alpha: A_K \to KG_K$ and $\beta: A_L \to L G_L$ are defined to be equivalent if and only if, for all finite-dimensional $kG$-supermodules $M$, the module $\alpha^*(M_K)$ is of finite flat dimension if and only if $\beta^*(M_L)$ is of finite flat dimension.

In \cite[\S9]{Benson:2022}, BIKP show for an arbitrary finite unipotent $k$-supergroup scheme $G$, there is a bijection $\Phi_G$ between $\Proj \Hbul(G,k)$ and the set $\Pi(G)$ of equivalence classes of $\pi$-points in $G$.\footnote{There are some inaccuracies at the beginning of \cite[\S9]{Benson:2022} in the discussion of the results from our paper \cite{Drupieski:2019b}, though these can be corrected by strategically replacing instances of the functor $\Hom_{\sGrp}(\Mr,G)$ with either $C_r(G)$ or $\Hom_{\sHopf}(\Pr,kG)$. The discussion is essentially correct as-is once the authors specialize to finite connected unipotent supergroups, which is their main case of interest. See Remark \ref{rem:notation-caution}, Remark \ref{rem:Hom(Mr,G)}, and \eqref{eq:Vrs(G)-stabilize} for details.} They also define the $\pi$-support of a module, and show that $\pi$-support identifies with support defined via cohomology.

\begin{question}
In light of the fact that the elementary supergroup schemes are inadequate for detecting projectivity for non-unipotent finite supergroup schemes: What modifications should be made, if any, to the definition of $\pi$-points for an arbitrary (non-unipotent) finite supergroup scheme, to enable an identification between $\pi$-support and cohomological support?
\end{question}
	
\makeatletter
\renewcommand*{\@biblabel}[1]{\hfill#1.}
\makeatother

\bibliographystyle{eprintamsplain}
\bibliography{survey-of-supports}

\end{document}